\documentclass[10pt, final, journal, letterpaper, oneside, twocolumn]{IEEEtran}

\makeatletter

\let\proof\@undefined
\let\endproof\@undefined

\let\theorem\@undefined
\let\endtheorem\@undefined
\makeatother

\usepackage{amssymb}
\usepackage{amsfonts}
\usepackage{latexsym, amssymb,amsmath, amsbsy,amsopn, amstext}
\usepackage{graphicx}
\usepackage{amsmath, amsbsy,bm}
\usepackage{hyperref}
\usepackage{cancel, color}
\usepackage{cite}
\usepackage[none]{hyphenat}
\usepackage{float}
\usepackage{ifpdf}
\usepackage{xcolor}
\usepackage{tikz}
\usetikzlibrary{arrows,shapes,snakes,shadows,positioning,automata,patterns}
\usetikzlibrary{trees,decorations.pathmorphing,decorations.markings}

\def\0{{\bf 0}}

\newtheorem{theorem}{Theorem}
\newtheorem{lemma}{Lemma}

\newtheorem{remark}{Remark}
\newtheorem{alg}{Algorithm}
\newtheorem{assumption}{Assumption}

\allowdisplaybreaks

\title{ Distributed  GNE seeking under partial-decision information over networks   via a doubly-augmented operator splitting approach 
}
\author{ Lacra Pavel
\thanks{This work was supported by NSERC Discovery Grant (261764).}
\thanks{Lacra Pavel is  with Department of Electrical and Computer Engineering, University of Toronto, Canada.
        {\tt \small pavel@control.utoronto.ca}}%
}

\begin{document}
\maketitle

\begin{abstract}
We consider distributed computation of  {\it generalized} Nash equilibrium (GNE) over networks, in games with shared coupling constraints. Existing methods require that each player has full access to opponents' decisions. In this paper, we assume that players have only partial-decision information, and can communicate with their neighbours over an arbitrary undirected graph.  We recast the problem as that of finding a zero of a sum of monotone operators through primal-dual analysis. To distribute the problem, we doubly augment variables, so that each player has  local decision estimates  and local copies of Lagrangian multipliers. We introduce a single-layer algorithm, fully distributed with respect to both primal and dual variables. We show its convergence to a variational GNE with fixed step-sizes, by reformulating it as a forward-backward  iteration for a pair of doubly-augmented monotone operators.  

\end{abstract}

\section{Introduction}


{\it Generalized Nash equilibrium} (GNE) problems in games with shared coupling constraints arise in various network scenarios where a set of players (agents) compete for limited network resources, e.g. power grids and smart grids,  \cite{Basar2012}, 
 optical networks, \cite{pavel2}, wireless communication networks, \cite{asuman,shanbhag4},   electric vehicle charging, \cite{lygeros2}. 
 The study of GNE dates back to  \cite{debreu,rosen}; a historical review is provided in \cite{faccinei1}, 
  \cite{pang}. Distributed GNE computation in monotone games has seen an increasing interest  in recent years,  \cite{shanbhag4,zhuminghui,payoff,lygeros2,grammatico_1,liangshu,grammatico_2,GrammaticoCDC2017,GrammaticoECC2018,YiPavelCDC1_2017,yipeng2TCNS}. 
Most works assume that each player has access 
 to all other agents' decisions - the classical setting of  {\em full-decision information}, either by observation or by a central node coordinator.  
 
There are many current networked applications where agents may only access or observe the decisions of their neighbours, and there is no central node to provide them with global information, i.e., a {\em partial-decision information} setting. The assumption of information exchange is motivated in networks where there is no central node that has bidirectional communications with all players to provide them with global information, as in peer-to-peer networks.  Application scenarios range from spectrum access in cognitive radio networks, where 
users adaptively adjust their operating parameters based on interactions with the environment and other users in the network, \cite{Cheng2014},  congestion games in ad-hoc networks, \cite{Tekin2012}, to networked Nash-Cournot competition, \cite{cournotgame}, and opinion dynamics in social networks, \cite{Srikant2014}, \cite{Ozdaglar2016}. 
 These examples are non-cooperative in the way decisions are made (each agent minimizes its own cost function), while agents exchange locally information with neighbours to compensate for the lack of global information on others' decisions. The first results on distributed NE seeking under such partial-decision information  have been for finite-action games, \cite{kar1}, and for aggregative games with {\em no coupling constraints}, \cite{shanbhag1}. Results were extended to  general continuous-kernel games in \cite{pavel4,GadjovPavelTAC2018,SalehiPavelIFAC_2017,ShiPavelACC_2017,YeHu_TAC_2017}, 
 for NE seeking problems only, in games with \emph{no coupling constraints}.  
 Inspired by work on NE seeking under partial-decision information, \cite{shanbhag1}, and by the recent elegant, operator-theoretic approach to GNE problems,  \cite{GrammaticoECC2018}, \cite{YiPavelCDC1_2017}, in this paper we consider GNE seeking in games with affine coupling constraints,  under partial-decision and local information exchange over an arbitrary  network. 
  

\emph{Literature review:} 
Distributed  (variational) GNE computation is an active research area, but  existing results are for 
 the classical setting of  {\em full-decision information}. 
Initial results were developed based on a variational inequality (VI) approach, \cite{faccinei1}, \cite{shanbhag4}.  For (pseudo)-monotone games, \cite{shanbhag4} adopts a single-layer Tikhonov regularization primal-dual algorithm, \cite{zhuminghui} proposes a primal-dual gradient approach,  while \cite{payoff} proposes a payoff-based algorithm,  all with diminishing step-sizes. Recently, an operator-splitting approach has proved to be very powerful;  it allows the design of GNE algorithms  that are guaranteed to globally converge with fixed step-sizes, with concise convergence proofs. 
Most results are for  {\em aggregative games}, \cite{lygeros2},  \cite{liangshu,grammatico_2,GrammaticoCDC2017,GrammaticoECC2018}. In \cite{grammatico_2,GrammaticoCDC2017,GrammaticoECC2018}, algorithms are semi-decentralized, requiring a central node (coordinator) to broadcast the common multipliers and/or aggregative variables, hence a star topology. 
This is relaxed in  \cite{liangshu} by  combining a  continuous-time consensus dynamics and a projected gradient,  still for aggregative games.  
For games with {\em generally coupled costs} and affine coupling constraints, distributed and center-free GNE seeking  is investigated via an operator approach  in \cite{YiPavelCDC1_2017,yipeng2TCNS,yipeng2}: a forward-backward algorithm, convergent  in strongly monotone games \cite{YiPavelCDC1_2017,yipeng2}, and  
preconditioned proximal algorithms for monotone games \cite{yipeng2TCNS}. Players communicate the local multipliers  over a network with arbitrary topology, in a distributed, peer-to-peer manner, but 
 each agent has access to the decisions of all other agents that influence his cost, hence  {\em full-decision information}. 

{\em Contributions:}  Motivated by the above, in this paper we consider distributed GNE seeking 
in a {\em partial-decision information} setting via an operator-splitting approach. 
We propose a fully distributed GNE seeking algorithm for games with generally coupled costs and affine coupling constraints, over networks with an arbitrary topology. To the best of our knowledge,  this is the first such algorithm in the literature.   
Based  on a primal-dual analysis of the variational inequality KKT conditions,  we reformulate the problem as that of finding zeros of  a sum of monotone operators and use the Laplacian matrix to distribute the computations. Different from \cite{YiPavelCDC1_2017} (perfect opponents' decision information), herein 
we distribute both the primal and the dual variables. 
To account for partial-decision information, we endow each agent with an auxiliary variable that estimates the other agents' decisions (primal variables), as in NE seeking over networks, \cite{GadjovPavelTAC2018,SalehiPavelIFAC_2017}.  
Compared to  \cite{YiPavelCDC1_2017,yipeng2},  this introduces technical challenges, as a change in an estimate induces a nonlinear change in an agent's dynamics. We make use of two selection matrices and we incorporate the Laplacian in an appropriate manner to do double duty, namely  to enforce  consensus of the local decision estimates (primal variables) and of  the local  multipliers (dual variables). Compared to  \cite{PavelCDC2018}, here we relax the assumption of cocoercivity of the extended pseudo-gradient. Under Lipschitz continuity  of the extended pseudo-gradient, we prove convergence 
 with  {\em fixed step-sizes} over any connected graph, by  leveraging monotone operator-splitting techniques, \cite{combettes1}. Specifically,  we reformulate the algorithm as a forward-backward iteration for doubly-augmented monotone operators, and   
 distribute the resolvent operation via a doubly-augmented metric matrix.  


The paper is organized as follows.
Section \ref{sec_notation} gives the notations and preliminary background.
Section \ref{sec_game_and_algorithm} formulates the game. Section \ref{sec_algorithm_develpo}   introduces the distributed GNE seeking algorithm and 
reformulates it as an operator-splitting  iteration. The convergence analysis is presented in  Section \ref{sec_alg_converge}, numerical simulations in  Section \ref{sec:numerical} and 
concluding remarks are given in  Section \ref{sec_concluding}. 
Some of the proofs are placed in the appendix.

\vspace{-0.225cm}
\section{Preliminary background}\label{sec_notation}

\vspace{-0.15cm}
\noindent {\it Notations. } For a vector $x  \!\in \! \bm{R}^m$, 
$x^T$  denotes its transpose 
and $\|x\|= \sqrt{x^Tx}$  the norm induced by inner product $\langle\cdot,\cdot\rangle$.
For a symmetric positive-definite matrix $\Phi$, $\Phi \succ 0$, $s_{\min}(\Phi)$ and $s_{\max}(\Phi)$ denote its minimum and maximum eigenvalues. The $\Phi$-induced inner product is  $\langle x,  y\rangle_{\Phi}=\langle \Phi x,y\rangle$ and the $\Phi$-induced norm, $\|x \|_\Phi=\sqrt{\langle \Phi x,x\rangle}$. 
For a matrix $A\in \bm{R}^{m\times n}$, let $\|A\| = \sigma_{\max} (A)$ denote  
the 2-induced matrix norm, where $\sigma_{\max} (A)$ is its maximum singular value.  
Let $\bm{1}_m\!=\!(1,...,1)^T  \!\in \! \bm{R}^m$ and
$\bm{0}_m \!=\!(0,...,0)^T \! \in \! \bm{R}^m$. For $\mathcal{N}\!=\! \{1,...,N\}$, $col(x_i)_{i\in  \mathcal{N}}$ or $[x_i]_{i\in  \mathcal{N}}$  denotes 
 the stacked vector obtained from vectors $x_i$, 
$diag((A_i)_{i\in  \mathcal{N}})$ 
the block diagonal matrix with  $A_1, . . . ,A_N$ on the
main diagonal. 
$Null(A)$ and $Range(A)$ are the null and range space of matrix $A$, respectively, while $[A]_{ij}$
stands for its $(i,j)$ entry. $I_m$ denotes the identity matrix in $\bm{R}^{m\times m}$.
 Denote 
 $\times_{i=1,...,N}\Omega_i$ or $\prod_{i=1}^N \Omega_i$ as the Cartesian product of the sets $\Omega_i,i=1,...,N$.

\vspace{-0.25cm}
\subsection{Monotone  operators}
The following are  from  \cite{combettes1}.  Let $\mathfrak{A} \! \!:\!\!\bm{R}^m \!\!\rightarrow \! \!2^{\bm{R}^m}$ be a set-valued operator.  
The domain of $\mathfrak{A}$ is $dom\mathfrak{A} \!= \! \{x\in \bm{R}^m| \mathfrak{A}x \neq \emptyset\}$ where $\emptyset$ is the empty set, and the range of $\mathfrak{A}$ is $ran\mathfrak{A}=\{y \in \bm{R}^m| \exists x, y\in \mathfrak{A}x\}$. The graph of $\mathfrak{A}$ is $gra\mathfrak{A}=\{(x,u) \in \bm{R}^m\times \bm{R}^m| u\in \mathfrak{A}x\}$;   the inverse of $\mathfrak{A}$ is defined through its graph as $gra\mathfrak{A}^{-1}=\{(u, x)| (x, u)\in gra \mathfrak{A}\}$.
The zero set of  $\mathfrak{A}$ is $zer\mathfrak{A}=\{x\in \bm{R}^m | \bm{0} \in \mathfrak{A}x\}$.
$\mathfrak{A}$ is called monotone if
$\forall (x,u), \forall(y,v)\in gra\mathfrak{A}$, 
$\langle x-y, u-v\rangle \geq 0.$
It is maximally monotone if $gra\mathfrak{A}$ is not {\it strictly} contained in the graph of any other monotone operator. The resolvent of $\mathfrak{A}$ is $\mathcal{J}_{\mathfrak{A}}=({\rm Id}+\mathfrak{A})^{-1}$, where ${\rm Id}$ is the identity operator. $\mathcal{J}_\mathfrak{A}$ is single-valued and $dom\mathcal{J}_{\mathfrak{A}}=\bm{R}^m$ if $\mathfrak{A}$ is maximally monotone. 
The composition of  $\mathfrak{A}$ and $\mathfrak{B}$ is denoted by $ \mathfrak{A}\circ \mathfrak{B}$. 
The sum $\mathfrak{A}+\mathfrak{B}$ is defined as $gra (\mathfrak{A}+\mathfrak{B})=\{(x,y+z)| (x,y)\in gra \mathfrak{A}, (x,z)\in gra \mathfrak{B}\}$.
If  $\mathfrak{A}$ and $\mathfrak{B}$ are maximally monotone operators and
$0\in int(dom \mathfrak{B}-dom \mathfrak{A})$, then $\mathfrak{A}+\mathfrak{B}$ is also maximally monotone. If $\mathfrak{B}$ is single-valued, then
$zer(\mathfrak{A}+\mathfrak{B})=Fix  (\mathcal{J}_{\mathfrak{A}}\circ({\rm Id}-\mathfrak{B}))$,  \cite[Prop.  25.1]{combettes1}, where $FixT$ denotes the set of fixed points of $T$. 

For a proper {\it lower semi-continuous convex} (l.s.c.) function $f$, its sub-differential $\partial f: domf\rightarrow 2^{\bm{R}^m}$ is $x \mapsto \{g\in \bm{R}^m| f(y)\geq f(x)+ \langle g, y -x \rangle, \forall y\in domf\}.$
$\partial f$ is a maximally monotone operator.
$Prox_{f}= \mathcal{J}_{\partial f}:\bm{R}^m\rightarrow dom f$, $ Prox_{f} :  x \mapsto  \arg\min_{u\in dom f } f(u)+\frac{1}{2} \|u -x \|^2$  is the proximal operator of $f$. Define the indicator function of  $\Omega$ as $\iota_{\Omega}(x)= 0$ if $ x\in \Omega$  and $\iota_{\Omega}(x)= \infty$ if $x\notin \Omega.$
For a closed convex set $\Omega$, $\iota_{\Omega}$ is a proper l.s.c. function and $\partial \iota_{\Omega}$ is the normal cone operator of  $\Omega$, 
   $N_{\Omega}(x)=\{v| \langle v, y -x\rangle\leq 0, \forall y\in \Omega\}  $. 

An operator $T: \Omega \subset \bm{R}^m\rightarrow \bm{R}^m$  is nonexpansive if
it is $1-$Lipschitz, i.e., $\|T(x)-T(y) \| \leq \|x-y \|, \forall x,y \in \Omega$. 
$T$ is $\alpha-$averaged ($\alpha\! \in \! (0,1)$),  if there exists a nonexpansive operator $T^{'}$ such that  $T\!=\!(1\!-\!\alpha){\rm Id}\!+\!\alpha T^{'}$.
By \cite[Prop. 4.25]{combettes1}, given 
$\alpha\in (0,1)$, $T\in \mathcal{A}(\alpha) $, where $\mathcal{A}(\alpha)$ denotes the class of $\alpha-$averaged operators, if and only if $\forall x, y \in \Omega$:\\
(i): $\|\!Tx\!-\!Ty \!\|^2 \!\leq \!\|x\!-\!y \|^2\! -\!\frac{1\!-\!\alpha}{\alpha} \| (x\!-\!y)\!-\!(Tx\!-\!Ty) \|^2$.\\
(ii): $\|Tx\!-\!Ty \|^2\!+\! (1\!-\!2\alpha)\|x\!-\!y\|^2 \leq 2(1\!-\!\alpha) \langle x\!-\!y, Tx\!-\!Ty \rangle$. \\ 
If $T \! \!\in \! \!\mathcal{A}(\frac{1}{2})$, $T$ is also called firmly nonexpansive. 
If  $\mathfrak{A}$ is maximally monotone, $\mathcal{J}_{\mathfrak{A}}\!\!=\!\!({\rm Id}\!+\!\mathfrak{A})^{-1}\!$ is firmly nonexpasive, \cite[\!Prop. \! 23.7]{combettes1}. 
Let the projection of $x$ onto $\Omega$ be $P_{\Omega}(x)\!=\!\arg\min_{y\in \Omega} \|x-y \|$, with 
$\!P_{\Omega}(x)\!=\!Prox_{\iota_{\Omega}}(x) \!=\!\mathcal{J}_{N_{\Omega}}(x)\!$. If $\Omega$ is closed and convex, $P_{\Omega}$  is firmly nonexpansive   since 
$N_{\Omega}$ is maximally monotone \cite[Prop.  4.8]{combettes1}.  
$T$ is called $\beta-$cocoercive if  $\beta T \in \mathcal{A}(\frac{1}{2})$, for $\beta >0$, i.e., 
$\beta\|T(x)-T(y) \|^2 \leq \langle x-y, T(x)-T(y) \rangle, \forall x, y \in \Omega$. 
If $f$  is convex differentiable, with $\theta-$Lipschitz gradient $\nabla f $, then $\nabla f $ is  $\frac{1}{\theta}-$cocoercive (cf. Baillon-Haddad theorem, \cite[Thm. 18.15]{combettes1}). 

\vspace{-0.25cm}
\subsection{Graph theory}

The following are from \cite{god}. Let  graph
$\mathcal{G}_c=(\mathcal{N},\mathcal{E})$ describe the information exchange among a set  
$\mathcal{N}$
of agents, where  $\mathcal{E} \! \subset \! \mathcal{N} \!\times \!\mathcal{N} $ is the edge set. 
If agent $i$ can get  information from agent $j$, then $(j,i) \in \mathcal{E}$ and
agent $j$ belongs to agent $i$'s  neighbour set $\mathcal{N}_i=\{ j | (j,i) \in
\mathcal{E}\}$,   $i \notin \mathcal{N}_i$.
$\mathcal{G}_c$ is undirected when
$(i,j)\in \mathcal{E}$ if and only if $(j,i)\in \mathcal{E}$.
$\mathcal{G}_c$ is connected if any two agents are
connected. Let  $W=[w_{ij}]\in \bm{R}^{N\times N}$  be the weighted adjacency matrix, with $w_{ij} >0$ if  
$j\in \mathcal{N}_i$ and $w_{ij}=0$ otherwise, and 
$Deg= diag ((d_i)_{i\in \mathcal{N}})$,  where $d_i=\sum_{j=1}^N w_{ij}$. Assume
$W=W^T$. 
The weighted Laplacian of $\mathcal{G}_c$ is
$L=Deg-W.$ 
When $\mathcal{G}_c$ is  connected and undirected, 0 is a simple
eigenvalue of  $L$, 
$L \bm{1}_N=\bm{0}_N$, $\bm{1}^T_{N} L=\bm{0}^T_N$; all other eigenvalues are positive.
Let the eigenvalues of $L$   in ascending order be $0<s_2(L)\leq ... \leq s_N(L)$, $d^* \leq s_N(L) \leq 2d^*$, where  $d^*=\max_i \{d_i\}$ is the maximal weighted degree. 

\section{Game formulation}\label{sec_game_and_algorithm}

Consider a group of agents (players) $\mathcal{N}=\{1,..., N \}$, where  
each player $i\in \mathcal{N}$ controls its local decision (strategy or decision) $x_i\in \bm{R}^{n_i}$.  
 Denote $x=col(x_i)_{i\in  \mathcal{N}} \in   \bm{R}^n$ as the 
decision profile,  i.e., the stacked vector of all the agents' decisions where $\sum_{i=1}^N n_i=n$. 
 We also write $x=col(x_i)_{i\in  \mathcal{N}}$ as $x=(x_i,x_{-i})$ where $x_{-i}=col(x_j)_{j\in  \mathcal{N}\setminus\{i}  = col(...,x_{i-1},x_{i+1},...)$
denotes the decision profile of all agents' decisions except player $i$. 
Agent $i$ aims to optimize its objective function $J_i(x_i,x_{-i})$,  
coupled to other players' decisions,  with respect to its own decision $x_i$ over its feasible decision set. 
Let the globally shared, affine coupled constrained set be
\begin{equation}\label{equ_coupling_set}
K := \prod_{i=1}^N \Omega_i \bigcap
\{ x \in \bm{R}^n | \sum_{i=1}^N A_ix_i \leq \sum_{i=1}^N b_i\}.
\end{equation}
where $\Omega_i \subset \bm{R}^{n_i}$ is a private feasible  set of player $i$, and $ A_i \in \bm{R}^{m\times n_i}$, $b_i\in \bm{R}^m$ its local data. 
Let  $\Omega=\prod_{i=1}^N \Omega_i$. 
A jointly-convex game with coupled constraints is represented by the set of inter-dependent optimization problems
\begin{equation}\label{GM}
\forall i \in \mathcal{N}, \, \min_{x_i\in \bm{R}^{n_i}} \;  J_i(x_i,x_{-i})   \;  \quad s.t. \;\; x_i\in  K_i(x_{-i}).
\end{equation}
where $K_i(x_{-i}):=\{ x_i\in  \bm{R}^{n_i}: (x_i,x_{-i}) \in K \}$ is the feasible decision set of agent $i$. 
A generalized Nash equilibrium (GNE) of game  \eqref{GM}, \eqref{equ_coupling_set}  is a profile $x^*=col(x^*_i)_{i\in  \mathcal{N}}$ 
at the intersection of all best-response sets,\vspace{-0.1cm}
\begin{equation}
x_i^* \in \arg\min J_i(x_i,x^*_{-i}) \; s.t. \; x_i\in  K_i(x^*_{-i}), \forall i\in \mathcal{N}.
\end{equation}
\begin{assumption}\label{assum1}
For each player $i$,  
$J_i(x_i,x_{-i})$ is continuously differentiable and convex in $x_i$, given  $x_{-i}$, and $\Omega_i$ is non-empty  compact and convex.  $K$ is non-empty and satisfies  Slater's constraint qualification.  
\end{assumption}
Denote $A=[A_1,...,A_N]$ and $b=\sum_{i=1}^N b_i$. Suppose $x^*$ is a GNE of game \eqref{GM}, \eqref{equ_coupling_set}  then for agent $i$, $x_i^*$ is the optimal solution to
the following convex optimization problem:
\begin{equation}\label{op_1}
\min_{x_i\in \bm{R}^{n_i}}  J_i(x_i,x^*_{-i}), \quad s.t.\; x_i\in \Omega_i, A_ix_i \leq b-\sum_{j \neq i,j\in \mathcal{N}} A_j x^*_j.\vspace{-0.1cm}
\end{equation}
A primal-dual characterization can be obtained via a Lagrangian for each agent $i$, 
\begin{equation}\label{local_lagrangian}
L_i(x_i,\lambda_i;x_{-i}) =  J_i(x_i,x_{-i})+\lambda^T_{i}(A x -b).
\end{equation}
with dual variable (multiplier) $\lambda_i\in \bm{R}^m_{+}$. When  $x_i^*$ is an optimal solution to \eqref{op_1}, there exists $\lambda_i^*\in \bm{R}_{+}^m$ such that the following  KKT conditions are satisfied:\vspace{-0.1cm}
\begin{equation}
\begin{array}{lll}\label{kkt_3}
\bm{0}_{n_i} =  \nabla_{x_i} L_i(x^*_i,\lambda^*_i;x^*_{-i}), \, x^*_i \in \Omega_i, \, \, i \in \mathcal{N} \\
\langle \lambda^*_i, Ax^* -b \rangle=0, \, -(Ax^* - b) \geq \bm{0}, \, \lambda_i^*\geq \bm{0},
\end{array}
\end{equation}\vspace{-0.1cm}
Equivalently, using the normal cone operator, \vspace{-0.1cm}
\begin{equation}
\begin{array}{lll}\label{kkt_1}
\bm{0}_{n_i} \in \nabla_{x_i} J_i(x^*_i,x^*_{-i}) + A_i^T \lambda_i^* +N_{\Omega_i}(x^*_i), i \in \mathcal{N}  \\
\bm{0}_m \in - (Ax^*-b) + N_{\bm{R}^m_{+}}(\lambda_i^*)
\end{array}
\end{equation}
Denote $\bm{\lambda}=col(\lambda_i)_{i\in  \mathcal{N}} $. 
By  \cite[Thm. 8, \S 4]{faccinei1} when $(x^*,\bm{\lambda}^*)$ satisfies KKT conditions \eqref{kkt_1},  $x^*$ is a GNE of game  \eqref{GM}, \eqref{equ_coupling_set}.

A GNE with the same Lagrangian multipliers for all the agents is called {\em variational GNE}, \cite{faccinei1}, 
which has the economic interpretation of no price discrimination,  \cite{shanbhag2}. A {\it variational  GNE}  of  game \eqref{GM}, \eqref{equ_coupling_set} is defined as  $x^* \in K$   solution of the following $VI(F,K)$:
\begin{equation}\label{vi}
 \langle F(x^*), x -x^*\rangle \geq 0, \, \,  \forall x\in K,
\end{equation}
where $F$ is  the {\em pseudo-gradient} of the game defined as:
\begin{equation}\label{pseudogradient}
F(x) = col(\nabla_{x_i} J_i(x_i,x_{-i}))_{i\in  \mathcal{N}}.
\end{equation}
$x^*$ solves $VI(F,K)$ if and only if there  exists a $\lambda^*\in \bm{R}^m$ such that the  KKT conditions are satisfied, \cite[\S 10.1]{FacchineiBOOK}, 
\begin{equation}
\begin{array}{lll}\label{kkt_2_BIG}
\bm{0}_n  &\in F(x^*) + A^T \lambda^* +N_{\Omega}(x^*) \\
\bm{0}_m & \in  -(Ax^*-b) + N_{\bm{R}^m_{+}}(\lambda^*) 
\end{array}
\end{equation}
where $N_{\Omega}(x^*)=\prod_{i=1}^N N_{\Omega_i}(x^*_i)$, or  component-wise, 
\begin{equation*}
\begin{array}{lll} 
\bm{0}_{n_i} \in \nabla_{x_i} J_i(x^*_i,x^*_{-i}) + A_i^T \lambda^* +N_{\Omega_i}(x^*_i),i\in \mathcal{N}  \\
\bm{0}_m \in  -(Ax^*-b) + N_{\bm{R}^m_{+}}(\lambda^*).  
\end{array}
\end{equation*}
Assumption \ref{assum1} 
 guarantees existence of a  solution to $VI(F,K)$  \eqref{vi}, by  \cite[Cor. 2.2.5]{FacchineiBOOK}. By  \cite[Thm. 9, \S 4]{faccinei1},  
every solution $x^*$ of  $VI(F,K)$  \eqref{vi} is a GNE of game \eqref{GM}. Furthermore, if $x^*$ together with $\lambda^*$ satisfies
 the KKT conditions  \eqref{kkt_2_BIG} for $VI(F,K)$  \eqref{vi},  then $x^*$  satisfies
the KKT conditions \eqref{kkt_1} with $\lambda_1^*\!=\!...\!=\!\lambda^*_N\!=\!\lambda^*$,  hence $x^*$ is variational GNE of game \eqref{GM}. 

Our aim is to design an iterative algorithm that finds a variational GNE under \emph{partial-decision information} over a network with arbitrary topology $\mathcal{G}_c$, by using an operator-theoretic approach. We first review typical iterative algorithms under full-decision information, where each agent has access to the others' decisions.

\subsection{Iterative Algorithm under Full-Decision Information}
\begin{assumption}\label{strgmon_Fassump}
 	$F$ is strongly monotone and Lipschitz continuous: there exists $\mu>0$ and $\theta_0 >0$ such that for any pair of points $x$ and $x'$,  
	$\langle x-x', F(x) - F(x') \rangle \geq \mu  \| x- x'\|^2$  and   
	$\| F(x) - F(x') \| \leq \theta_0  \| x- x'\|$.	
\end{assumption} 
Strong monotonicity  of $F$ is a standard assumption under which convergence of projected-gradient type algorithms is guaranteed with fixed step-sizes, e.g.  \cite{zhuminghui},\cite{lygeros2}, \cite{YiPavelCDC1_2017}, \cite{GrammaticoECC2018}. 
Under Assumption  \ref{assum1}, \ref{strgmon_Fassump}, the $VI(F,K)$, \eqref{vi}, has  a unique solution $x^*$ (cf. \cite[Thm. 2.3.3]{FacchineiBOOK}),  thus the game \eqref{GM} has a unique variational GNE. Assuming each player has access to the others' decisions $x_{-i}$, i.e., {\it full-decision information}, a primal-dual projected-gradient GNE algorithm  is \vspace{-0.2cm}  
\begin{equation}\label{dal_1_semi_decentralized}
\begin{array}{ll}
x_{i,k+1}   & =  P_{\Omega_i}\big( x_{i,k}-\tau ( \nabla_{x_i} J_i(x_{i,k},x_{-i,k}) + A_i^T \lambda_{k} )\big) \\
\lambda_{k+1} & = P_{\bm{R}^{m}_{+}}\Big (\lambda_{k} + \sigma\big ( A (2x_{k+1} -x_{k}) - b \big)\Big), \vspace{-0.4cm}
\end{array}
\end{equation}
where $x_{i,k}$,  $\lambda_{k}$ denote $x_i$, $\lambda$ at iteration $k$ and $\tau$ and $\sigma$ are fixed step-sizes. The dual variable $\lambda$ is handled by a center (coordinator) as in \cite{GrammaticoECC2018} hence \eqref{dal_1_semi_decentralized} is  semi-decentralized. 

Algorithm \eqref{dal_1_semi_decentralized}  is an instance of an operator-splitting method for finding zeros of a sum of monotone operators, \cite[\S 25]{combettes1}. To see this, note that  the KKT conditions \eqref{kkt_2_BIG} can be written as $col(x^*, \lambda^*) \in zer  \mathcal{T}$ where the operator $\mathcal{T}$ is defined by the concatenated right-hand side of   \eqref{kkt_2_BIG}. $\mathcal{T} $ can be  
split  as $\mathfrak{A} + \mathfrak{B}$, where 
  operators $\mathfrak{A}$, $\mathfrak{B}$  are defined as \vspace{-0.2cm}
\begin{align}\label{operator_tilde_A_B}
& \mathfrak{A}:\; 
\left [
                                     \begin{array}{c}
                                         x \\
                                          \lambda
                                        \end{array}
                                      \right]
\mapsto
\left [
                                     \begin{array}{c}
                                         N_{\Omega}(x) \\
                                         N_{\bm{R}^m_{+}}(\lambda)
                                        \end{array}
                                      \right] +
\left [
                                     \begin{array}{cc}
                                         0 & A^T \\
                                         - A & 0
                                        \end{array}
                                      \right]
\left [
                                     \begin{array}{c}
                                         x \\
                                          \lambda
                                        \end{array}
                                      \right] \nonumber
\\                                    
& \mathfrak{B}:\; 
\left [
                                     \begin{array}{c}
                                         x \\
                                          \lambda
                                        \end{array}
                                      \right]
\mapsto
\left [
                                     \begin{array}{c}
                                        F(x) \\
                                          b
                                        \end{array}
                                      \right]
\end{align}
Algorithm \eqref{dal_1_semi_decentralized} can be obtained as a forward-backward iteration, \cite[\S 25.3]{combettes1}, for zeros of $\Phi^{-1}\mathfrak{A} + \Phi^{-1} \mathfrak{B}$, where  $\Phi = \left [
                                     \begin{array}{cc}
                                         \tau^{-1} & -A^T \\
                                         - A & \sigma^{-1}
                                        \end{array}
                                      \right]
$ is a metric matrix. We note that different GNE seeking algorithms can be obtained for different splitting of $\mathcal{T}$, 
with convergence conditions dependent on monotonicity properties of  $\mathfrak{A}$ and $\mathfrak{B}$. Notice that $\mathfrak{A}$ \eqref{operator_tilde_A_B} is  maximally monotone (similar arguments for this can be found in Lemma \ref{lem_monotone}), and 
under Assumption \ref{strgmon_Fassump}, $\mathfrak{B}$ is cocoercive. Convergence of \eqref{dal_1_semi_decentralized} to $x^*$, $\lambda^*$ can be proved  for sufficient conditions on the fixed-step sizes such that $\Phi \succ 0$. 

\section{Distributed Algorithm under Partial-Decision Information}\label{sec_algorithm_develpo}  %

In this section we consider a  {\em partial-decision information} setting, where the agents do not have full information on the others' decisions $x_{-i}$. 
 We propose an  algorithm that 
  allows agents to 
 find 
a  variational GNE based on local information exchange with neighbours, over a communication graph $\mathcal{G}_c$ with arbitrary topology, under the following assumption.
\vspace{0.2cm}

\begin{assumption}\label{connectivity}
$\mathcal{G}_c$ is undirected and connected. 
\end{assumption}

\vspace{0.2cm}

{\it Our approach is based on the interpretation of the KKT conditions \eqref{kkt_2_BIG} as a zero-finding problem of a sum of operators. To deal with (incomplete) partial-decision information and to distribute the computations, we introduce estimates and lift the original problem to a higher-dimensional space. This space (called the augmented space), is  doubly-augmented (in both primal  and dual variables), and the original space is its consensus subspace.   
We appropriately define a pair of doubly-augmented operators,  such that any zero 
of their sum lies on the consensus subspace, and has variational GNE $x^*$ and $\lambda^*$ as its components. }

\vspace{0.2cm}
We describe next the algorithm variables.  
Agent $i$ controls its local decision $x_i\in \bm{R}^{n_i}$, and a local copy of multiplier (dual variable) $\lambda_i \in \bm{R}_{+}^m$ for the estimation of $\lambda^*$  in \eqref{kkt_2_BIG}. To cope with partial-decision information, we endow each player  with an auxiliary variable $\mathbf{x}^i $ that provides an estimate of other agents' primal variables (decisions), as done in \cite{GadjovPavelTAC2018,SalehiPavelIFAC_2017} for NE seeking. Thus agent $i$ maintains  $\mathbf{x}^i =col(\mathbf{x}^i_j)_{j \in \mathcal{N}} \in \bm{R}^n$,  
where $\mathbf{x}^i_j $ is player $i$'s estimate of player $j$'s decision and $\mathbf{x}^i_i = x_i$ is its decision. Note that $\mathbf{x}^i = (x_i, \mathbf{x}^i_{-i})$, where $\mathbf{x}^i_{-i}$ represents player $i$'s estimate vector without its own decision $x_i$. In steady-state all estimates should be equal, i.e., $\mathbf{x}^i=\mathbf{x}^j$ and $\lambda_i=\lambda_j$.  \textit{Each agent uses  the relative feedback from its neighbours such that  in steady-state these estimates, on both {\it primal} and {\it dual variables}, agree one with another.}  An additional local auxiliary variable $z_i\in \bm{R}^m$ is used for the  coordination needed to  satisfy the  coupling constraint and to reach  consensus of the  local multipliers (dual variables) $\lambda_i$. 
Agents  exchange local $\{ \mathbf{x}^i, \lambda_i,z_i\}$ via the arbitrary topology communication graph $\mathcal{G}_c=\{\mathcal{N},\mathcal{E}\}$. $(j,i)\in \mathcal{E}$ if player $i$ can receive $\{ \mathbf{x}^j, \lambda_j,z_j\}$ from player $j \in \mathcal{N}_i$, where $\mathcal{N}_i=\{j| (j,i)\in \mathcal{E}\}$  denotes its set of neighbours.



\vspace{0.3cm}

The distributed algorithm for player  $i$ is given as follows. 

\vspace{0.2cm}
\begin{alg}\label{dal_1}
\\
\noindent\rule{0.49\textwidth}{0.7mm}
Initialize: $x_{i,0}\!\in \Omega_i$,  $\!\mathbf{x}^i_{-i,0}\! \in \!\bm{R}^{n-n_i}$, $\lambda_{i,0}\!\in\! \bm{R}_{+}^m$, $\!z_{i,0} \!\in \!\bm{R}^m$.\\
Iteration: \vspace{-0.4cm}
\begin{equation}
\begin{array}{ll}
x_{i,k+1} \!\! \!\!   & =  P_{\Omega_i}\left (x_{i,k}-\tau_i \big ( \nabla_{x_i} J_i(x_{i,k},\mathbf{x}^i_{-i,k}) +A_i^T \lambda_{i,k} \right .  \\
\!\! \!\!&\qquad \qquad \qquad  \quad \left .+{c} \sum_{j\in \mathcal{N}_i} w_{ij} (x_{i,k}- \mathbf{x}^j_{i,k}) \big )\right )  \\
\mathbf{x}^i_{-i,k+1} \!\! \!\!   & =  \mathbf{x}^i_{-i,k}  - \tau_i {c} \sum_{j\in \mathcal{N}_i} w_{ij}  \big (\mathbf{x}^i_{-i,k}- \mathbf{x}^j_{-i,k} \big ) \\
z_{i,k+1} \!\!  \!\!  & =  z_{i,k}                + \nu_i  \sum_{j\in \mathcal{N}_i} w_{ij}(\lambda_{i,k}-\lambda_{j,k}) \\
\lambda_{i,k+1}\!\! \!\!  & = P_{\bm{R}^{m}_{+}}\Big(\lambda_{i,k} + \sigma_i\big( A_i(2x_{i,k+1} -x_{i,k})-b_i
\nonumber\\
\!\! \!\!&\; -\sum_{j\in \mathcal{N}_i} w_{ij}(2(z_{i,k+1}-z_{j,k+1})-(z_{i,k}-z_{j,k}))
\nonumber\\
\!\! \!\!&\;\quad\quad  \qquad \qquad \quad  -\sum_{j\in\mathcal{N}_i} w_{ij}(\lambda_{i,k}-\lambda_{j,k})\big)\Big)\label{al_3}
\end{array}
\end{equation}
\noindent\rule{0.49\textwidth}{0.7mm}
Here $x_{i,k}$, $\!\mathbf{x}^i_{-i,k} $, $\!z_{i,k}$, $\!\lambda_{i,k}$ denote $x_i$,  $\!\mathbf{x}^i_{-i}$, $\!z_i$, $\lambda_i$ at iteration $k$,
$c\!>\!0$ is a design parameter, $\tau_i,\!\nu_i,\! \sigma_i \!>\!0$ are fixed step-sizes of player $i$, and $W\!=[w_{ij}]$ the weighted adjacency matrix of $\mathcal{G}_c$.
\end{alg}

\begin{remark}
The update for $x_i$ in Algorithm \ref{dal_1} employs a projected-gradient descent of the local Lagrangian function in \eqref{local_lagrangian} with an extra proportional term  of  the consensual errors (disagreement) between his primal variables and his neighbours' estimates. 
The updates for $ \mathbf{x}^i_{-i}$  and  $z_i$  can be regarded as discrete-time integrations for the consensual errors of local decision estimates and dual variables. 
Finally, $\lambda_i$ is updated  by a  combination of the projected-gradient ascent of local Lagrangian \eqref{local_lagrangian} and a proportional-integral term for consensual errors. 
 Each player  knows only its local data in game \eqref{GM},   $J_i$, $\Omega_i$, $A_i$ and $b_i$,  own private information, i.e., cost function, preference and decision ability. 
  $A_i$ characterizes how  agent $i$ is involved in the coupled constraint (shares the global resource), assumed to be privately known by player $i$. The globally shared constraint $A x\leq b$ couples the agents' feasible decision sets, but is {\it not known} by any agent.
\end{remark}

\begin{remark}\label{rem_Alg1_double_aug}
Compared to  algorithm \eqref{dal_1_semi_decentralized}, Algorithm \ref{dal_1} is completely distributed  (without any central coordinator), i.e., primal-distributed and dual-distributed  over $\mathcal{G}_c$. The algorithms in  \cite{YiPavelCDC1_2017},  \cite{GadjovPavelTAC2018} are special cases of Algorithm \ref{dal_1}. When each agent has  access to all players' decisions that affect its cost, the estimates 
$ \mathbf{x}^i_{-i}$  are not needed (set $c\!=\!0$), and Algorithm \ref{dal_1}  reduces to the dual-distributed, perfect-information case one in  \cite{YiPavelCDC1_2017} (dual distributed).  On the other hand, in a game with no coupling constraints (set $A\!=\!0$, $b\! =\!0$), the $\lambda_i$ (hence the $z_i$) are not needed, and  Algorithm \ref{dal_1} reduces to a discrete-time version of the primal-distributed dynamics in \cite{GadjovPavelTAC2018}.  


\end{remark}

Next, we write Algorithm \ref{dal_1} in compact form, using two matrices to manipulate the selection of agent $i$'s decision variables, $x_i$, and estimate variables, $\mathbf{x}^i_{-i}$. Let \vspace{-0.2cm}
\begin{align} \label{eq:actualStratREMatrix}
\mathcal{R}_i & =  \begin{bmatrix}
\mathbf{0}_{n_i\times n_{<i}} && & I_{n_i} & && \mathbf{0}_{n_i\times n_{>i}}
 \end{bmatrix} \\
\mathcal{S}_i &= \begin{bmatrix}
	I_{n<i} & \mathbf{0}_{n<i \times n_i} & \mathbf{0}_{n<i \times n>i} \\
	\mathbf{0}_{n>i \times n<i} & \mathbf{0}_{n>i \times n_i} & I_{n>i} 
\end{bmatrix}\end{align}
where   $n_{<i} = \sum_{j<i\ j\in\mathcal{N}} n_j$, $n_{>i} = \sum_{j>i\ j\in\mathcal{N}} n_j$. Hence  $\mathcal{R}_i$ selects the $i$-th $n_i$-dimensional component  from an $n$-dimensional vector, 
while $\mathcal{S}_i$ removes it.   Thus,  $\mathcal{R}_i \, \mathbf{x}^i \!=\!\mathbf{x}^i_i\!=\!x_i$ and $\mathcal{S}_i \, \mathbf{x}^i\!=\!\mathbf{x}^i_{-i}$. 
 With $x_i=\mathcal{R}_i \, \mathbf{x}^i$, the stacked decisions $x\!=\!col(x_i)_{i \in \mathcal{N}} \!\in \! \mathbf{R}^n$ 
can be written as $x \!= \!\mathcal{R}\mathbf{x} $, where   $\mathcal{R} \! =\! diag((\mathcal{R}_i)_{i \in \mathcal{N}})$ and $\mathbf{x} \!=col(\mathbf{x}^i)_{i \in \mathcal{N}} \in \mathbf{R}^{Nn}$. 
 Similarly,  the stacked estimates are $ col(\mathbf{x}^i_{-i})_{i \in \mathcal{N}} = \mathcal{S}\mathbf{x} \in  \mathbf{R}^{(N-1)n}$, where  
   $\mathcal{S}  = diag((\mathcal{S}_i)_{i \in \mathcal{N}})$.  
These two matrices, $\mathcal{R}$ and $\mathcal{S}$, play a key role in the following.  Using \eqref{eq:actualStratREMatrix}, it can be seen that  both  are full row rank and moreover,  \vspace{-0.2cm}
\begin{align}\label{eq:IdRSMatrix} 
\mathcal{R}^T\mathcal{R} + \mathcal{S}^T\mathcal{S} = I_{Nn}, \quad \quad \nonumber  \\
\,\,	\mathcal{R}\mathcal{S}^T = \mathbf{0}_{n}, \, \, \mathcal{S}\mathcal{R}^T = \mathbf{0}_{Nn-n}  \\ 	
\,\,	\mathcal{R}\mathcal{R}^T = I_{n}, \, \, \mathcal{S}\mathcal{S}^T = I_{Nn-n}.\nonumber
\end{align}
Furthermore, with   $\mathbf{v} \!:=\mathcal{S}\mathbf{x}$, we can write $\mathbf{x}= \mathcal{R}^T x + \mathcal{S}^T \mathbf{v}$.   
With these notations, we write Algorithm \ref{dal_1} in stacked form, using boldface notation for stacked  variables (all local copies).



\begin{lemma}\label{lem_ALG_compact}
Let $x_k\!\!= \!\!col(x_{i,k})_{i \in \mathcal{N}}$, 
$\mathbf{x}_k\!\!=\!\! col(\mathbf{x}^i_{k})_{i \in \mathcal{N}}$, 
$\bm{\lambda}_k\!\!=\!\!col(\lambda_{i,k})_{i \in \mathcal{N}}$, 
$\bm{z}_k\!\!=\!\! col(z_{i,k})_{i \in \mathcal{N}}$. 
Then, Algorithm \ref{dal_1} is equivalently written in stacked notation as \vspace{-0.2cm}
\begin{align}
x_{k+1}   &=  P_{\Omega}\big(x_{k}-\bm{\tau}_x  (  \mathbf{F}(\mathbf{x}_{k}) +\bm{A}^T \bm{\lambda}_{k}+ {c} \mathcal{R}\mathbf{L}_x \,  \mathbf{x}_k)\big) \label{Alg1BLOCK_1} \\\mathcal{S}\mathbf{x}_{k+1}   &= \mathcal{S}\mathbf{x}_{k}- \bm{\tau}_s \, {c} \mathcal{S}\mathbf{L}_x \, \mathbf{x}_k \label{Alg1BLOCK_2}\\
\bm{z}_{k+1}   &=  \bm{z}_{k}                + \bm{\nu}  \mathbf{L}_{\lambda} \, \bm{\lambda}_{k} \label{Alg1BLOCK_3}\\
\bm{\lambda}_{k+1} &= P_{\bm{R}^{Nm}_{+}}\big(\bm{\lambda}_{k} + \bm{\sigma} ( \bm{A} (2x_{k+1} -x_{k})-\bm{b} \label{Alg1BLOCK_4} \\
\quad \, \, & \qquad \qquad \qquad - \mathbf{L}_{\lambda} \, \bm{\lambda}_{k}  -\mathbf{L}_{\lambda} \, (2\bm{z}_{k+1}-\bm{z}_k))\big)\nonumber 
\end{align}
where $ \mathbf{F}$  is defined as \vspace{-0.2cm}
\begin{equation}\label{pseudogradient_EKT}
\mathbf{F}(\mathbf{x})=  col(\nabla_{x_i} J_i(x_i,\mathbf{x}^i_{-i}))_{i \in \mathcal{N}}, 
\end{equation}  
$\bm{A}\! \! = \!diag(\!(A_i)_{i \in \mathcal{N}}\!$, 
$\!\!\bm{b}\!\! \! = \!col(b_i)_{i \in \mathcal{N}}\!$, $\!\mathbf{L}_{x} \!\! \! =\!L \!\otimes \!I_n$, $\!\mathbf{L}_{\lambda} \! \!\! =L \!\otimes \!I_m$, 
$\mathcal{R}  \!\!= \!diag(\!(\mathcal{R}_i)_{i \in \mathcal{N}}\!)$, 
$\mathcal{S}  = diag(\!(\mathcal{S}_i)_{i \in \mathcal{N}}\!)$, 
$\bm{\tau}_x   \!\!=  \!\!diag\!(\tau_i I_{n_i})_{i \in \mathcal{N}}\!$, 
$\!\bm{\tau}_s \!\!  =  \!\!diag\!(\tau_i I_{n-n_i})_{i \in \mathcal{N}}\!)$, 
$\!\bm{\nu}  \!\!  = \!\! diag\!(\nu_i I_{m})_{i \in \mathcal{N}}\!$, 
 $\!\bm{\sigma}  \!\!  = \!\! diag\!(\sigma_i I_{m})_{i \in \mathcal{N}}\!$. 
\end{lemma}\hfill $\Box$
\begin{remark}
In Algorithm \ref{dal_1},  instead of evaluating its gradient at actual decisions, as in  $\nabla_{x_i} J_i(x_{i},x_{-i})$,  
each player  evaluates its  gradient at local estimates, $\nabla_{x_i} J_i(x_i, \mathbf{x}^i_{-i})$. 
The stacked form  $ \mathbf{F}(\mathbf{x})$, \eqref{pseudogradient_EKT}, called the {\em extended pseudo-gradient}, is the extension of $F$, \eqref{pseudogradient} to the augmented space of decisions and estimates. 
When 
  these estimates are identical, $\mathbf{x}^i = x$ for all $i$, then $ \mathbf{F}(\mathbf{1}_N \otimes x) = F(x)$. 
\end{remark}



Based on Lemma \ref{lem_ALG_compact}, we show  next that  
Algorithm \ref{dal_1} can be written as a forward-backward iteration for finding zeros of the sum of two doubly-augmented operators $\Phi^{-1} {\bm{\mathfrak{A}}}$ and $\Phi^{-1} {\bm{\mathfrak{B}}}$, where ${\bm{\mathfrak{A}}}$, ${\bm{\mathfrak{B}}}$  are related to $\mathfrak{A}$, $\mathfrak{B}$,  \eqref{operator_tilde_A_B}, and $\Phi$ is a (preconditioning) metric matrix.   
Let $\varpi\!=\! col(\mathbf{x}, \bm{z},\bm{\lambda}) \! \in \! {\bm{\Omega}}$, where ${\bm{\Omega}}\!:=\! \bm{R}^{Nn} \! \times \! \bm{R}^{Nm}\! \times  \!\bm{R}_{+}^{Nm}$.  
%
Define  ${\bm{\mathfrak{A}}}: {\bm{\Omega}}\rightarrow 2^{\bm{R}^{Nn+2Nm}}$,  $ {\bm{\mathfrak{B}}}: {\bm{\Omega}}\rightarrow \bm{R}^{Nn+2Nm}$ as  
\begin{align}\label{op_hat_A}
& {\bm{\mathfrak{A}}}: 
\varpi \mapsto
\left [\!\!\!
  \begin{array}{c}
   \mathcal{R}^T N_{\Omega}(\mathcal{R}\mathbf{x}) \!\!\\  
    \bm{0} \\
 N_{\bm{R}^{Nm}_{+}}(\bm{\lambda})      		\!\!\!\!
  \end{array}
\! \right] \! + \!
\left[\!
  \begin{array}{ccc}
\!\!\!\!    \bm{0} & \bm{0} & \mathcal{R}^T \bm{A}^T\!\! \\
\!\!\!\!    \bm{0} & \bm{0} & -\mathbf{L} _\lambda \!\! \\
 \!\!\!\!   - \bm{A} \mathcal{R} & \mathbf{L}_\lambda       & \bm{0} \!\! \\
  \end{array}
\! \right] \! \!
\varpi
 \nonumber
\\
& {\bm{\mathfrak{B}}}: 
\varpi \mapsto  
 \left [\!
  \begin{array}{c}
   \mathcal{R}^T  \mathbf{F}(\mathbf{x}) +  {c} \mathbf{L}_x \mathbf{x} \\  
     \bm{0} \\
   \mathbf{L}_\lambda \, \bm{\lambda} +\bm{b} 	
    \end{array}
\! \right]
\end{align}
where $N_{\Omega}(\mathcal{R}\mathbf{x}) =N_{\Omega}(x)= \prod_{i=1}^N N_{\Omega_i}(x_i)$,  
 $\mathcal{R}^T N_{\Omega}(\mathcal{R}\mathbf{x})=\{ \mathcal{R}^T v \, |  v \in N_{\Omega}(\mathcal{R}\mathbf{x})\} $, 
 $N_{\bm{R}^{Nm}_{+}}(\bm{\lambda})=\prod_{i=1}^N N_{\bm{R}^m_{+}}(\lambda_i)$. 
 
 Let  the matrix $\Phi$ be defined as
\begin{equation}\label{metric_matrixBis}
{\Phi}=\left[\!
  \begin{array}{ccc}
    \bm{\tau}^{-1} & \bm{0} & -\mathcal{R}^T \bm{A}^T \\
    \bm{0} & \bm{\nu}^{-1} & \mathbf{L}_\lambda  \\
    -\bm{A} \mathcal{R} & \mathbf{L}_\lambda &  \bm{\sigma}^{-1} \\
  \end{array}
\!\right].
\end{equation}
where  $\!\bm{\tau} \!=\! diag((\tau_i I_n)_{i\in \mathcal{N}})$,  $\!\bm{\tau}^{-1} \!=\! diag((\tau^{-1}_i I_n)_{i\in \mathcal{N}})$, and  
$\bm{\nu}^{-1}\!, \bm{\sigma}^{-1}$ are similarly defined from $\bm{\nu}, \bm{\sigma}$. 

\begin{lemma}\label{lem_fix_approximation}
Let  $\varpi_k=col({\mathbf{x}}_k, {\bm{z}}_k,{\bm{\lambda}}_k)$, $\bm{\mathfrak{A}}$, $\bm{\mathfrak{B}}$,  $\Phi$ as in \eqref{op_hat_A},  \eqref{metric_matrixBis}. Suppose that ${\Phi} \!\succ \!0$ and  ${\Phi}^{-1} \!{\bm{\mathfrak{A}}}$  is maximally monotone. Then the following hold:\\
(i):  Algorithm \ref{dal_1} is equivalent to \vspace{-0.15cm}
\begin{equation} \label{compact_operator_1}
- {\bm{\mathfrak{B}}}(\varpi_{k}) \in {\bm{\mathfrak{A}}}(\varpi_{k+1})+{\Phi}(\varpi_{k+1}-\varpi_{k}), 
\end{equation}
\vspace{-0.6cm} 
\begin{equation}\label{com_fix} 
\hspace{-0.2cm}\text{or,} \, \, \varpi_{k+1}\!=\! ({\rm Id}\!+\!{\Phi}^{-1} {\bm{\mathfrak{A}}})^{-1}\!\!\circ\! ({\rm Id}\!-\! {\Phi}^{-1}{\bm{\mathfrak{B}}}) \varpi_{k}\!:=\! T_2 \!\circ \! T_1 \varpi_k,
\end{equation} 
where $T_1:={\rm Id}-{\Phi}^{-1}{\bm{\mathfrak{B}}} $ and $T_2:=({\rm Id}+{\Phi}^{-1}{\bm{\mathfrak{A}}})^{-1}$. \\
(ii):  Any limit point $\overline{\varpi}=col(\overline{\mathbf{x}}, \overline{\bm{z}},\overline{\bm{\lambda}})$  of Algorithm \ref{dal_1} is 
  a zero of ${\bm{\mathfrak{A}}}+{\bm{\mathfrak{B}}}$ and a fixed point of $T_2 \circ T_1$. 
\end{lemma}\hfill $\Box$

\begin{remark}\label{rem_FB_AB_B_Phi}
Algorithm 1,  written as   \eqref{com_fix}  is a {\em forward-backward iteration} for finding zeros of ${\Phi}^{-1} {\bm{\mathfrak{A}}} \!+ \! {\Phi}^{-1} {\bm{\mathfrak{B}}}$, or   fixed-point iteration for  $T_2 \!\circ \! T_1$ 
 \cite[\S 25.3]{combettes1}.  
It alternates a forward step $({\rm Id}-{\Phi}^{-1}{\bm{\mathfrak{B}}})$, 
and a backward step  $({\rm Id}+{\Phi}^{-1}{\bm{\mathfrak{A}}})^{-1}$. 
 Typically,  the backward step evaluates the resolvent of a monotone operator, while the forward step evaluates a cocoercive operator. Note that we put the skew-symmetric part  (monotone but not cocoercive) in ${\bm{\mathfrak{A}}}$,  \eqref{op_hat_A}, and $\mathbf{F}$, $\mathbf{L}_x$ and $  \mathbf{L}_\lambda$  in ${\bm{\mathfrak{B}}}$.   This separation between constraint data and game cost functions data will be instrumental later on to study monotonicity and  cocoercivity properties of the two operators. 
Note that with the standard choice $\Phi = I$,  the resolvent  $({\rm Id}+{\bm{\mathfrak{A}}})^{-1}$ needs evaluated,  which cannot be done in a distributed manner.  
This is where a non-identity metric matrix ${\Phi}$ as in \eqref{metric_matrixBis}  helps. Specifically, for ${\bm{\mathfrak{A}}}$  \eqref{op_hat_A}, let  ${\bm{\mathfrak{A}}}={\bm{\mathfrak{A}}}_1\!+\!{\bm{\mathfrak{A}}}_2$, where ${\bm{\mathfrak{A}}}_1=\mathcal{R}^TN_{\Omega} \mathcal{R}(\mathbf{x})\times \bm{0}_{Nm} \times N_{\bm{R}^{Nm}_{+}}(\bm{\lambda})$
 and ${\bm{\mathfrak{A}}}_2$ is the skew-symmetric matrix in \eqref{op_hat_A}. Then, with   ${\Phi}$ as in \eqref{metric_matrixBis},  \eqref{compact_operator_1} is written as 
$- {\bm{\mathfrak{B}}}(\varpi_{k}) \in {\bm{\mathfrak{A}}}_1(\varpi_{k+1})+ ({\bm{\mathfrak{A}}}_2  +\Phi)\varpi_{k+1}-\Phi \varpi_{k}
$, 
where 
$
{\bm{\mathfrak{A}}}_2 \!+\! {\Phi} \!=\! \left[\!\!
  \begin{array}{ccc}
    \bm{\tau}^{-1} & \!\bm{0} &\! \bm{0} \\
    \bm{0} & \!\bm{\nu}^{-1} & \!\bm{0}  \\
    -2\bm{A} \mathcal{R} & \!2\mathbf{L}_\lambda & \! \bm{\sigma}^{-1} \\
  \end{array}
\!\!\right]
$. Since this is lower block-triangular, 
 the resolvent can be evaluated distributively via only projection, matrix multiplying and local communication. Using $P_{\bm{R}_{+}^{Nm}}(\bm{\lambda})=({\rm Id}\!+\!N_{\bm{R}_{+}^{Nm}})^{-1}$, it can be checked that the explicit iterations are as in \eqref{Alg1BLOCK_1}-\eqref{Alg1BLOCK_4}. 
\end{remark}
 \begin{remark}\label{rem_op_hat_A_B}
Operators   $\bm{\mathfrak{A}}$, $\bm{\mathfrak{B}}$, \eqref{op_hat_A}, are doubly-augmented extensions  of $\mathfrak{A}$, $\mathfrak{B}$,  \eqref{operator_tilde_A_B}, obtained by introducing local copies of primal and  dual variables, from $x  \in \bm{R}^n$ to $\mathbf{x} =[\mathbf{x}^i]_{i\in\mathcal{N}} \in \bm{R}^{Nn}$,  $\lambda \in \bm{R}^m$ to $\bm{\lambda} = [\lambda_i]_{i\in\mathcal{N}} \in \bm{R}^{Nm}$ 
and with auxiliary variables $\bm{z} =  [z_i]_{i\in\mathcal{N}} \in \bm{R}^{Nm}$. 
The Laplacian matrices  $\bm{L}_x$ and $\bm{L}_\lambda$   help to enforce the consensus  of the local primal variables $\mathbf{x}^i$ and of the  local  dual multipliers $\lambda_i$, as well as the feasibility of the affine coupling constraints. 
The auxiliary variables $z_i$ help for the consensus of $\lambda_i$s and to decouple the constraint, in the sense of estimating the contribution of the other agents in the constraint. 
 Note that if we set $\mathcal{R} \!=\!I$, $\mathbf{L}_x \!=\! 0$, $\mathbf{x} \!=\! \bm{1}_N \!\otimes \!x$   then $ \mathbf{F}(\mathbf{1}_N \! \otimes \! x) \!=\! F(x)$ and ${\bm{\mathfrak{A}}}$ and ${\bm{\mathfrak{B}}}$ collapse to the operators ${\overline{\mathfrak{A}}}$ and ${\overline{\mathfrak{B}}}$ in \cite{YiPavelCDC1_2017}, distributed in dual variables only. Furthermore, if we also set $\mathbf{L}_\lambda =0$ and take $\bm{\lambda} = \bm{1}_N \! \otimes \!\lambda$, matrix $\bm{A}$ reduces to $A$, and we recover  operators $\mathfrak{A}$  and $\mathfrak{B}$, \eqref{operator_tilde_A_B}. 
Extension to the partial-decision information case introduces technical challenges in the convergence analysis: coupling in $\mathbf{F}$ (augmented space), non-square $\mathcal{R}$, and the extra term involving  $\mathbf{L}_x$. We will exploit properties of these extra terms to derive properties for  the doubly-extended operators ${\bm{\mathfrak{A}}}$ and ${\bm{\mathfrak{B}}}$. 
 \end{remark}

\section{Convergence analysis}\label{sec_alg_converge}

In this section we prove the convergence of Algorithm \ref{dal_1}.   
First, in Theorem \ref{thm_zeroHAT_is_correct}, based on the fact that any limiting  point of Algorithm \ref{dal_1} is a zero of  ${\bm{\mathfrak{A}}}+{\bm{\mathfrak{B}}}$ (Lemma \ref{lem_fix_approximation}(ii)), we  characterize the zeros of ${\bm{\mathfrak{A}}}+{\bm{\mathfrak{B}}}$. We show that any zero is on the consensus subspace and solves the $VI(F,K)$  \eqref{vi}, thus any limiting  point of Algorithm \ref{dal_1} finds a variational GNE of game \eqref{GM}. 
{Let  $\mathbf{E}_x \! \! \!= \!\{  \!\mathbf{x}\! \in \!\mathbf{R}^{Nn} | \mathbf{x}^i \!\!=\! \mathbf{x}^j, \forall i,j \!\!\in\! \mathcal{N} \!\}\! \!= \!\{  
 \mathbf{x} \!=\! \mathbf{1}_N  \!\otimes  x,  
x \!\in \! \bm{R}^{n}\! \}$ denote the estimate consensus subspace,  $\mathbf{E}_x^\perp$ its orthogonal complement  with $\mathbf{R}^{Nn} \! \!=\! \mathbf{E}_x \oplus \mathbf{E}_x^\perp$,   
$\mathbf{E}_\lambda \! \!\!=\! \{ \! \bm{\lambda} \! \in \! \mathbf{R}^{Nm} \! |
\bm{\lambda} \!=\! \mathbf{1}_N \!\otimes \!\lambda,  \lambda \!\in \! \bm{R}^{m}\! \}$ the multiplier consensus subspace. Consider $\mathbf{L}_{x} \!\! \! =\!L \!\otimes \!I_n$, $\mathbf{L}_{\lambda} \! \!\! =L \!\otimes \!I_m$.  Under Assumption \ref{connectivity} on  $\mathcal{G}_c$, $Null (\mathbf{L}_x) \! = \! Range \{ \mathbf{1}_N  \!\otimes \!I_n \}\! \!=\!\mathbf{E}_x$,  $Range (\mathbf{L}_x) \!=\! Null (\mathbf{L}_x)^\perp  \!\!= \!  Null \{ \mathbf{1}^T_N \! \otimes \! I_n \}\!=\!\mathbf{E}_x^\perp$,   $Null (\mathbf{L}_\lambda) \!=\! Range \{ \mathbf{1}_N  \!\otimes \! I_m \}\!=\! \mathbf{E}_\lambda$, $Range (\mathbf{L}_\lambda)  \!\!=\! Null \{ \mathbf{1}^T_N \!\otimes \! I_m \}\!$.   

\begin{theorem}\label{thm_zeroHAT_is_correct}
Suppose that Assumptions \ref{assum1}-\ref{connectivity} hold. Consider operators $\bm{\mathfrak{A}}$,$\bm{\mathfrak{B}}$,  \eqref{op_hat_A}. Then the following statements hold. \\
(i): Given any $\overline{\varpi}^*\!:=\!col(\mathbf{x}^*,\!\bm{z}^*\!,\!\bm{\lambda}^*) \!\in \! zer({\bm{\mathfrak{A}}}\!+\!{\bm{\mathfrak{B}}})$, 
then   {$\mathbf{x}^* \! \in \! \mathbf{E}_x$  and  $\bm{\lambda}^*\!  \in \! \mathbf{E}_\lambda$}, with $\mathbf{x}^*\! = \!\bm{1}_N \! \otimes x^*$ and $\bm{\lambda}^*\!=\!\bm{1}_N \!\otimes  \lambda^*$, where $x^*$, $\lambda^*$ satisfy the  KKT conditions \eqref{kkt_2_BIG}, $col(x^*,\lambda^*)\!\in \! zer(\mathfrak{A}\!+\!\mathfrak{B})$, for  $\mathfrak{A}$, $\mathfrak{B}$, \eqref{operator_tilde_A_B}. Moreover, $x^*$ solves the  $VI(F,K)$  \eqref{vi}, hence $x^*$ is a variational GNE of game \eqref{GM}. 
(ii):  $zer({\bm{{\mathfrak{A}}}}+{\bm{{\mathfrak{B}}}}) \! \neq \emptyset$. 
\end{theorem}

{\bf Proof:} 
(i) Let  $\overline{\varpi}^*:=col(\!\mathbf{x}^*,\!\bm{z}^*,\!\bm{\lambda}^*\!) \!\in\! zer({\bm{\mathfrak{A}}}\!+\!{\bm{\mathfrak{B}}})$. By \eqref{op_hat_A}, \vspace{-0.45cm}
\begin{align}\label{ZERO_op_hat_AB}
&  \bm{0}_{Nn} \in  \mathcal{R}^T  \mathbf{F}(\mathbf{x}^*)  +  {c}\mathbf{L}_{x} \mathbf{x}^* 
 +\mathcal{R}^T \bm{A}^T   \bm{\lambda}^* + \mathcal{R}^T N_{\Omega}(\mathcal{R} \mathbf{x}^*) \nonumber \\
 &    \bm{0}_{Nm} = -\mathbf{L}_\lambda  \bm{\lambda}^* \\
 &    \bm{0}_{Nm} \in  \mathbf{L}_{\lambda}  \bm{\lambda}^*  + \bm{b} +N_{\bm{R}^{Nm}_{+}}(\bm{\lambda}^*)	 
      -\bm{A} \mathcal{R} \mathbf{x}^* +  \mathbf{L}_\lambda  \bm{z}^*  \nonumber
\end{align}
From the  {\em first} line it follows that for some $v \in  N_{\Omega}( \mathcal{R}\mathbf{x}^*)$   \vspace{-0.15cm} 
 \begin{align}\label{ZERO_op_hat_ABEQ}
&  \bm{0}_{Nn} =  \mathcal{R}^T \big [ \mathbf{F}(\mathbf{x}^*)  + \bm{A}^T   \bm{\lambda}^* +   v \big ]  +  {c}\mathbf{L}_{x} \mathbf{x}^*  
\end{align}
From the {\em second} line of \eqref{ZERO_op_hat_AB} it follows that  $\bm{\lambda}^* \in Null(\mathbf{L}_{\lambda})\!=\!\mathbf{E}_\lambda$, 
hence  $\bm{\lambda}^*  = \mathbf{1}_N \otimes \lambda^*$, for some $ \lambda^* \in \bm{R}^{m}$. Note  that $(\mathbf{1}_N^T\otimes I_n)\mathcal{R}^T  = I_n$ (by \eqref{eq:actualStratREMatrix}) and  $(\mathbf{1}_N^T\otimes I_n)\mathbf{L}_x = \mathbf{0}_{n\times Nn}$. 
The, premultiplying \eqref{ZERO_op_hat_ABEQ} by $(\mathbf{1}_N^T\otimes I_n)$ yields \vspace{-0.2cm}
\begin{align}\label{ZERO_op_hat_ABis2EQ}
&  \bm{0}_{n} =   \mathbf{F}(\mathbf{x}^*)   + \bm{A}^T   \bm{\lambda}^* +  v 
\end{align}
Substituting \eqref{ZERO_op_hat_ABis2EQ} into  \eqref{ZERO_op_hat_ABEQ} yields 
$ \bm{0}_{Nn}\! = \mathbf{L}_{x} \mathbf{x}^* $. Thus,   $\mathbf{x}^* \! \in \!Null(\mathbf{L}_{x})\! = \! \mathbf{E}_x$ 
(by Assumption \ref{connectivity}), 
 and  $\mathbf{x}^*  \!= \!\mathbf{1}_N \! \otimes \! x^*$, for some $ x^* \!\in \!\bm{R}^{n}$.  
Using $\!\mathbf{x}^* \! =\! \mathbf{1}_N \!\otimes \!x^*$ and $\bm{\lambda}^* \! =\! \mathbf{1}_N \!\otimes \!\lambda^*$ in \eqref{ZERO_op_hat_ABis2EQ} yields 
$ \bm{0}_{n} \!= \!   \mathbf{F}(\mathbf{1}_N \!\otimes \!x^*)  \!+ \!\bm{A}^T(\mathbf{1}_N \! \otimes \!\lambda^*) \!+ \!  v $ 
for  $v \in  N_{\Omega}( \mathcal{R}\mathbf{x}^*)$, or 
$  \bm{0}_{n} \in    \mathbf{F}(\mathbf{1}_N \otimes x^*)  + \bm{A}^T ( \mathbf{1}_N \otimes \lambda^* )+  N_{\Omega}(x^*)$. 
With $ \mathbf{F}(
\mathbf{1}_N \otimes x^*)  = F( x^*)$, 
  $ \bm{A}^T ( \mathbf{1}_N \otimes \lambda^* )= A^T  \lambda^*$, this is equivalent to 
$ \bm{0}_{n} \in   F(x^*)  + A^T  \lambda^* +  N_{\Omega}(x^*) $,
which is the \textit{first } line in \eqref{kkt_2_BIG}. 

Using $\mathbf{x}^*  \!= \! \mathbf{1}_N \! \otimes  \!x^*$, $\bm{\lambda}^* \! = \! \mathbf{1}_N \!\otimes \!\lambda^*$, $\mathbf{L}_\lambda  \bm{\lambda}^* \!= \!\bm{0}_{Nm} $ in the {\em third} line of \eqref{ZERO_op_hat_AB},   
yields 
$   \bm{0}_{Nm} \in 
  \bm{b} +N_{\bm{R}^{Nm}_{+}}(\mathbf{1}_N \otimes \lambda^*)	
      -\bm{A}  \mathcal{R} (\mathbf{1}_N \otimes x^*) +  \mathbf{L}_\lambda  \bm{z}^*. 
$  
Thus, with  $  \mathcal{R} ( \mathbf{1}_N \!\otimes  \! x^*) \! =\! x^*$,
\vspace{-0.25cm}
\begin{align}\label{4thZERO_op_hat_ABis3EQ} 
 &    \bm{0}_{Nm} = 
  \bm{b} + col (w_i)_{i \in \mathcal{N}} -     \bm{A}  x^* +  \mathbf{L}_\lambda  \bm{z}^* .
\end{align}
for some $w_i \! \in \! N_{\bm{R}^{m}_{+}}( \lambda^*)$, $i \!\in \!\mathcal{N}$. Premultiplying \eqref{4thZERO_op_hat_ABis3EQ}  by $(\mathbf{1}_N^T \!\otimes \! I_m)$, with $(\mathbf{1}_N^T \! \otimes \! I_m)\mathbf{L}_\lambda \!=\! \mathbf{0}_{m \!\times \!Nm}$,  $(\mathbf{1}_N^T \!\otimes  \!I_m)\bm{A}  \! =\! A$, yields that 
$
   \bm{0}_{m} =   \sum_{i =1}^N b_i + \sum_{i =1}^N w_i	  -      A x^*, 
$ or,  
$\bm{0}_{m}   \in    b  \!+\! \sum_{i =1}^N N_{\bm{R}^{m}_{+}}( \lambda^*)	\!-  \!    A x^*.$ 
This gives the \textit{second} line in \eqref{kkt_2_BIG},  using $N_{\bigcap_{i=1}^N \mathcal{S}_i}\!=\!\sum_{i=1}^N \! \! N_{\mathcal{S}_i}$ if $\bigcap_{i=1}^N int(\mathcal{S}_i)\!\!\neq \! \emptyset$, \cite[Cor. 16.39]{combettes1}.  
Thus,  the KKT  conditions \eqref{kkt_2_BIG} for $VI(F,K)$  \eqref{vi} are satisfied for $x^*$ and $\lambda^*$, $col(x^*,\lambda^*) \in zer(\mathfrak{\mathfrak{A}}+\mathfrak{\mathfrak{B}})$, hence 
$x^*$  is a variational GNE of game \eqref{GM}. 

(ii) Under Assumption \ref{assum1} and \ref{strgmon_Fassump}, $VI(F,K)$ \eqref{vi} has a unique solution $x^*$,  \cite[Thm. 2.3.3]{FacchineiBOOK}, hence   there exists  $\lambda^* \!\in \! \bm{R}^m$ such that  \eqref{kkt_2_BIG} holds, \cite[Prop.  1.2.1]{FacchineiBOOK},  $col(x^*,\lambda^*) \in zer(\mathfrak{A}+\mathfrak{B})$. 
The first two lines in  \eqref{ZERO_op_hat_AB} are satisfied with  $\mathbf{x}^* \! = \!\mathbf{1}_N \!\otimes \! x^*$ and $\bm{\lambda}^*\!=\!\bm{1}_N \! \otimes \!\lambda^*$. Using \eqref{kkt_2_BIG}  it can be shown that  $\exists \bm{z}^* \!\in \! \bm{R}^{Nm}$ such that the third line in \eqref{ZERO_op_hat_AB} is also satisfied, hence   $zer({\bm{{\mathfrak{A}}}}\!+\!{\bm{{\mathfrak{B}}}}) \!\neq \!\emptyset$. This is omitted due to space constraints (see the proof of  \cite[Thm. 1]{yipeng2} for similar arguments). 
\hfill $\Box$

Note that \eqref{ZERO_op_hat_AB} can be regarded as the  KKT conditions  \eqref{kkt_2_BIG},  doubly-augmented in both the primal and the dual space using local variables, $\mathbf{x}$, $\bm{\lambda}$, and  auxiliary variables $\bm{z}$. 

Next we focus on proving the convergence of Algorithm \ref{dal_1}, based on its interpretation as a forward-backward iteration for zeros of ${\Phi}^{-1}{\bm{\mathfrak{A}}} +{\Phi}^{-1}{\bm{\mathfrak{B}}}$,  or a  fixed-point iteration of $T_2 \!\circ \! T_1$  (see \eqref{com_fix}, Lemma \ref{lem_fix_approximation}). 
 This is done  in Theorem \ref{thm_convergence}, using nonexpansiveness properties of $T_1$ and $T_2$, obtained in Lemma \ref{lem_PhiA_B_prop} from properties of ${\Phi}^{-1}{\bm{\mathfrak{A}}}$ and ${\Phi}^{-1}{\bm{\mathfrak{B}}}$.  These are based on monotonicity properties of  $\bm{\mathfrak{A}}$ and $\bm{\mathfrak{B}}$ (Lemma \ref{lem_monotone}),  under a sufficient condition for $\Phi$ to be positive definite (Lemma \ref{lem_monotone_metric}). To prove Lemma \ref{lem_PhiA_B_prop} and \ref{lem_monotone},  a key result is Lemma \ref{boldRF_plus_cL_strgmon_lemma}, which shows  how a \textit{restricted} monotonicity property can be achieved in the augmented space under Lipschitz continuity of  $\mathbf{F}$. 
 \begin{assumption}\label{Lipchitz_boldFassump}
 	The extended pseudo-gradient  $\mathbf{F}$,  \! \eqref{pseudogradient_EKT},  is Lipschitz continuous: there exists $\theta>0$ such that for any $\mathbf{x}$ and $\mathbf{x'}$, 
$\| \mathbf{F}(\mathbf{x})-\mathbf{F}(\mathbf{x'}))\| \leq \theta \| \mathbf{x}- \mathbf{x'}\|$.    
	\end{assumption} 
\begin{remark}
 Note that Assumption~\ref{Lipchitz_boldFassump} on $\mathbf{F}$ is weaker than those in 
  \cite{SalehiPavelIFAC_2017,ShiPavelACC_2017,YeHu_TAC_2017}, \cite{PavelCDC2018}. 
In the classical, full-decision information setting,  convergence proof relies on cocoercivity  (strong monotonicity) of $F$, e.g. \cite{YiPavelCDC1_2017}.  Cocoercivity of $\mathbf{F}$ (the extension of cocoercivity of $F$ to the augmented space), is sometimes used in works on NE computation in partial-decision information setting,    \cite{SalehiPavelIFAC_2017,ShiPavelACC_2017,YeHu_TAC_2017}, \cite{PavelCDC2018}. 
Unlike  distributed optimization, where  cocoercivity  is automatically satisfied, in a game context, monotonicity of $\mathbf{F}$ is \emph{not automatically satisfied} in the augmented space, because of coupling  to the others' actions and because of partial convexity. It does hold  in games with shared coupling constraints but separable costs. We show next how, in generally coupled games, under Assumption~\ref{Lipchitz_boldFassump}, a \textit{restricted} monotonicity property can be achieved in the augmented space. 
\end{remark}

{\begin{lemma}\label{boldRF_plus_cL_strgmon_lemma} 
		Consider that Assumptions~\ref{assum1}-\ref{Lipchitz_boldFassump} hold and let  \vspace{-0.15cm} 
		\begin{equation}\label{eqboldRF_plus_cL_strgmon_5}
\Psi =  \left [
\begin{array}{ccc}
\frac{\mu}{N} && - \frac{\theta + \theta_0}{2\sqrt{N}}  \\ 		
- \frac{\theta + \theta_0}{2\sqrt{N}}  && c  \, s_2(L) - \theta 
\end{array}
\right ]
\end{equation}
Then,  for any $c > c_{\min}$, where $c_{\min} \, s_2(L) =  \frac{(\theta+\theta_0)^2}{4\mu} +\theta$,  $\Psi \succ 0$ 	and, for  any $ \mathbf{x}$  
and any  $ \mathbf{x'} \in \mathbf{E}_x$, \vspace{-0.2cm} 
\begin{eqnarray}\label{eqboldRF_plus_cL_strgmon}
&&\hspace{-1.5cm}(\!\mathbf{x}\!-\! \mathbf{x'})^T \! \big (\mathcal{R}^T(\mathbf{F}(\mathbf{x})\!-\!\mathbf{F}(\mathbf{x'}))\!+\!c \mathbf{L}_x(\mathbf{x} \!-\! \mathbf{x'}) \big )\!\geq \!\bar{\mu} \|\mathbf{x}\!-\!\mathbf{x'}\|^2,
\end{eqnarray}
where $\bar{\mu}\!:=\!s_{\min}(\Psi)>0$ and $\mathbf{E}_x\!=\!Null(\mathbf{L}_x)$ is the $n$- dimensional consensus subspace.  
\end{lemma}}\vspace{-0.3cm}\hfill$\Box$
\begin{remark}\label{rem_theta_quadratic} Lemma \ref{boldRF_plus_cL_strgmon_lemma} is instrumental to show global convergence of Algorithm \ref{dal_1} under Assumption \ref{Lipchitz_boldFassump}. 
In  Section \ref{sec:numerical} we consider a quadratic game example and check analytically Assumption  \ref{Lipchitz_boldFassump}. 
 \end{remark}}
 {Property \eqref{eqboldRF_plus_cL_strgmon} in Lemma \ref{boldRF_plus_cL_strgmon_lemma} means that, for a sufficiently large $c\!>\!0$, $ \mathcal{R}^T\mathbf{F}\!+ c\mathbf{L}_x$ is  $\bar{\mu}$-strongly monotone in a \textit{restricted} set of directions (since $\mathbf{x'} \in\mathbf{E}_x $), which we call $\bar{\mu}$-\textit{restricted} strongly monotone.   
This is a weaker monotonicity property,  
similar to the notion of \textit{restricted} convexity used in DOP, \cite{Zhang2017},  and high-dimensional statistical estimation,  \cite{NegahbanNIPS2009}.  
To show this key property, an instrumental step is the decomposition of the augmented  space $\mathbf{R}^{Nn}\! =\!\mathbf{E}_x \!\oplus \!\mathbf{E}_x^\perp$ into  the consensus subspace $\mathbf{E}_x\!=\!Null(\mathbf{L}_x)$ (where  $F$ is strongly monotone), and its orthogonal complement $\mathbf{E}^\perp\!=\!Null(\mathbf{L}_x)^\perp$ (where $ \mathbf{L}_x$ is strongly monotone), (see proof in Appendix). 
Based on this property, we show in Lemma \ref{lem_monotone} that the operator $\bm{\mathfrak{B}}$ is \textit{restricted} cocoercive, which  is sufficient to prove global convergence of Algorithm~\ref{dal_1}.} 
\begin{lemma}\label{lem_monotone}
Suppose Assumptions \ref{assum1}-{\ref{Lipchitz_boldFassump}} hold and {$c$ is selected such that $  c  > c_{\min}$, with $c_{\min}$ as in Lemma   \ref{boldRF_plus_cL_strgmon_lemma}. } 
Then the following hold for operators ${\bm{\mathfrak{A}}}$ and ${\bm{\mathfrak{B}}}$, \eqref{op_hat_A}. \\(i): ${\bm{\mathfrak{A}}}$ is maximally monotone. \\
(ii): ${\bm{\mathfrak{B}}}$   is $\beta$-\textit{restricted} cocoercive: for any $\varpi $ and any {$\varpi' \!\in \bm{\Omega}_E$,  where 
${\bm{\Omega}}_E\!=\! \mathbf{E}_x \!\times \!\! \bm{R}^{Nm}\! \! \times  \!\bm{R}_{+}^{Nm} $,}
 the following holds, \vspace{-0.2cm}\begin{equation}\label{eq:B_cocercive_restricted}
  \!\langle \varpi -\varpi',\bm{\mathfrak{B}} \varpi  - \bm{\mathfrak{B}} \varpi' \rangle  \geq  \beta  \|  \bm{\mathfrak{B}}\varpi- \bm{\mathfrak{B}}\varpi' \|^2  
\end{equation}
{where $\beta \in (0, \min \{ \frac{\bar{\mu}}{\bar{\theta}^2},  \frac{1}{2d^*}  \}]$,  $\bar{\theta} =  \theta + 2 c d^*$, $\bar{\mu}=s_{\min}(\Psi)$, $\Psi$ as  in \eqref{eqboldRF_plus_cL_strgmon_5} and $d^*$ is the maximal weighted degree of  $\mathcal{G}_c$.}
\end{lemma}\hfill $\Box$

The next result  follows from  Lemma 3 in \cite{ YiPavelCDC1_2017}) and  Gershgorin's theorem. 
\begin{lemma}\label{lem_monotone_metric} Consider ${\Phi}$, \eqref{metric_matrixBis}. Given any $\delta >0$,  if each player choose step-sizes $\tau_i, \nu_i, \sigma_i$ in Algorithm \ref{dal_1} such that  \vspace{-0.15cm} 
\begin{equation}
\begin{array}{lll}\label{step_size_choice}
0<&\tau_i   & \leq  ({\max_{j=1,...n_i} \{\sum^m_{k=1} | [A_i^T]_{jk} |\}+\delta})^{-1},\\
0<&\sigma_i & \leq ({\max_{j=1,...m} \{ \sum^{n_i}_{k=1} | [A_i]_{jk} | \}+2d_i+ \delta})^{-1}\\
0<&\nu_i    & \leq  ({2d_i+\delta})^{-1}
\end{array}
\end{equation}
then ${\Phi} \succ 0$ and  $ {\Phi} - \delta I_{n+2Nm} \succeq 0$. 
\end{lemma}\vspace{-0.1cm}\hfill $\Box$

The next result gives  properties of ${\Phi}^{-1}{\bm{\mathfrak{A}}}$ and $ {\Phi}^{-1}{\bm{\mathfrak{B}}}$ 
under the ${\Phi}-$induced norm $ \| \!\!\cdot \!\!\|_{{\Phi}}$. 

\begin{lemma}\label{lem_PhiA_B_prop}Suppose Assumptions \ref{assum1}-\ref{Lipchitz_boldFassump} hold. Take any {$c$ such that $  c  > c_{\min}$, with $c_{\min}$ as in Lemma   \ref{boldRF_plus_cL_strgmon_lemma}, } and any $ \delta > \frac{1}{2\beta}$, 
where $\beta \in (0, \min \{ \frac{\bar{\mu}}{\bar{\theta}^2},  \frac{1}{2d^*}  \}]$,  $\bar{\theta} =  \theta +2  c d^*$, $\bar{\mu}=s_{\min}(\Psi)$, $\Psi$ as  in \eqref{eqboldRF_plus_cL_strgmon_5} and $d^*$ is the maximal weighted degree of  $\mathcal{G}_c$. 
Suppose that the step-sizes $\tau_i, \nu_i, \sigma_i$ in Algorithm \ref{dal_1} are chosen to satisfy \eqref{step_size_choice}. 
Then, the following hold for the  operators ${\Phi}^{-1}{\bm{\mathfrak{A}}}$, ${\Phi}^{-1}{\bm{\mathfrak{B}}}$ and  $T_1={\rm Id}-{\Phi}^{-1}{\bm{\mathfrak{B}}},   T_2=({\rm Id}+{\Phi}^{-1}{\bm{\mathfrak{A}}})^{-1} $, under the $\Phi$-induced norm $\| \!\cdot \!\|_\Phi$, with ${\Phi}$ as in \!\eqref{metric_matrixBis}.   \\
(i) ${\Phi}^{-1}{\bm{\mathfrak{A}}}$ is  maximally monotone and  $T_2 \in \mathcal{A}(\frac{1}{2})$. \\
(ii)  ${\Phi}^{-1}{\bm{\mathfrak{B}}}$ is {$\beta\delta$-\textit{restricted} cocoercive and $T_1$  is \textit{restricted} nonexpansive, such that }
 {for any  $\varpi$ and any $\varpi' \in {\bm{\Omega}}_E$,  \vspace{-0.15cm} 
\begin{equation}\label{eq:prop1_T1}
\begin{array}{ll}
\| T_1 \varpi \!-\!  T_1 \varpi' \|^2_{\Phi} \!\!&\leq \!  \| \varpi - \varpi' \|^2_{\Phi}  \\
& \hspace{-1.3cm} - (2 \beta \delta -1) \| \varpi  - \varpi ' - \left (T_1 \varpi -  T_1 \varpi'\right) \|^2_{\Phi}
\end{array}
\end{equation}}
\end{lemma}\hfill $\Box$
The next result shows the  convergence of  Algorithm \ref{dal_1} based on its compact reformulation \eqref{com_fix} and properties of $T_1$, $T_2$. 
\begin{theorem}\label{thm_convergence}
Suppose Assumptions \ref{assum1}-\ref{Lipchitz_boldFassump} hold. 
Take any {$c$ such that $  c  > c_{\min}$, with $c_{\min}$ as in Lemma   \ref{boldRF_plus_cL_strgmon_lemma}, } and  any $ \delta > \frac{1}{2\beta}$, where $\beta \in (0, \min \{ \frac{\bar{\mu}}{\bar{\theta}^2},  \frac{1}{2d^*}  \}]$,  $\bar{\theta} =  \theta + 2 c d^*$, $\bar{\mu}=s_{\min}(\Psi)$, $\Psi$ as  in \eqref{eqboldRF_plus_cL_strgmon_5} and $d^*$ is the maximal weighted degree of  $\mathcal{G}_c$.  
If step-sizes $\tau_i, \nu_i, \sigma_i$ in Algorithm \ref{dal_1}  are chosen satisfy \eqref{step_size_choice}, 
then, 
for each player $i$ 
its $\{\mathbf{x}^i_k\}_{k \in \mathbf{N}}$ converges to the same value 
$x^*$  variational GNE  of game \eqref{GM}, and 
 its   local decision $x_{i,k} = \mathbf{x}^i_{i,k}$ converges 
to the corresponding component in  $x^*$, while the local  multiplier $\lambda_{i,k}$ of all agents converge to the same multiplier corresponding to the KKT condition \eqref{kkt_2_BIG}, i.e.,
$\lim_{k\rightarrow \infty} x_{i,k}= x_i^*,  \quad \lim_{k\rightarrow \infty}  \lambda_{i,k} =\lambda^*$, $\forall i \in \mathcal{N}$. 
\end{theorem}
{\bf Proof:} 

By  Lemma \ref{lem_monotone_metric}, $\Phi \succ 0$, and by Lemma \ref{lem_PhiA_B_prop}, ${\Phi}^{-1}{\bm{\mathfrak{A}}}$ is maximally monotone and ${\Phi}^{-1}{\bm{\mathfrak{B}}}$ is \textit{restricted} $\beta\delta-$cocoercive. By Lemma \ref{lem_fix_approximation}(ii), Algorithm \ref{dal_1} can be written in the compact form  \eqref{com_fix}, fixed-point iteration $\varpi_{k+1}=T_2 \circ T_1\varpi_{k}$. 
Consider {
any  $\varpi^* \in zer({\bm{\mathfrak{A}}}+{\bm{\mathfrak{B}}})$}, or equivalently any fixed point of $T_2 \circ T_1$, i.e., $\varpi^* =T_2 T_1 \varpi^*$. By Lemma \ref{lem_fix_approximation} and \eqref{com_fix} it follows that  \vspace{-0.2cm}
\begin{equation}
\begin{array}{lll}\label{equ_thm_5_2}
\|\varpi_{k+1}- \varpi^*\|_{\Phi} \!\!&= \|  T_2( T_1 \varpi_{k})-  T_2 (T_1 \varpi^*)\|_{\Phi} \\
& \hspace{-1cm} \leq  \| T_1 \varpi_{k}-   T_1 \varpi^* \| \leq  \| \varpi_{k}-  \varpi^* \|_{\Phi}
\end{array}
\end{equation}
where  the first inequality follows  from $T_2$ being nonexpansive (by  Lemma  \ref{lem_PhiA_B_prop}(i)), and the second one from \eqref{eq:prop1_T1}, Lemma  \ref{lem_PhiA_B_prop}(ii), since $\varpi^* \!\in\! \bm{\Omega}_E$ (cf. Theorem \ref{thm_zeroHAT_is_correct}(i)) and $2\beta\delta \! >\!1$ by assumption. 
 Hence the sequence $\{\|\varpi_{k}-\varpi^* \|_{\Phi}\}$ is  non-increasing and bounded from below.
 By the monotonic convergence theorem and since $\Phi \succ 0$,  $\{\| \varpi_{k}-  \varpi^* \| \}$ is bounded  and converges for every $\varpi^* \in zer({\bm{\mathfrak{A}}}+{\bm{\mathfrak{B}}})$.\\
Denote $\xi = \frac{1}{2\delta \beta} \in (0,1)$.  Using again \eqref{com_fix}, we can write \vspace{-0.2cm}
\begin{equation}
\begin{array}{lll}\label{equ_thm_57_1}
&&\hspace{-0.8cm}\|\varpi_{k+1}- \varpi^* \|_{\Phi}^2= \|  T_2 ( T_1 \varpi_{k})-  T_2  (T_1 \varpi^*) \|_{\Phi}^2\\
&& \leq  \|T_1 \varpi_{k}- T_1 \varpi^*\|_{\Phi}^2 \\
&&\quad - \| (  T_1 \varpi_{k}-  T_1 \varpi^*) - (T_2 T_1 \varpi_{k}-  T_2  T_1 \varpi^*)\|^2_{\Phi} \\
&& \leq \|\varpi_{k}-  \varpi^*  \|^2_{\Phi} \\
&& \quad - \| (  T_1 \varpi_{k}-  T_1 \varpi^*) - (T_2 T_1 \varpi_{k}-  T_2  T_1 \varpi^*)\|^2_{\Phi} \\
&& \quad -\frac{1-\xi}{\xi}\| \varpi_{k}-\varpi^*- (T_1 \varpi_{k}- T_1 \varpi^*) \|^2_{\Phi}
\end{array}
\end{equation}
where the first inequality follows by \cite[Prop. 4.25(i)]{combettes1} for $T_2\in \mathcal{A}(\frac{1}{2})$, and the second  one from \eqref{eq:prop1_T1}, Lemma  \ref{lem_PhiA_B_prop}(ii), for $T_1$ with  $2\delta \beta\! =1/\xi$. Next, we use,  (cf.  \cite[Cor. 2.14]{combettes1}), \vspace{-0.2cm}
\begin{equation}\label{equation_norm}
\alpha \|x\|^2 + (1-\alpha)\|y \|^2 = \|\alpha x +(1-\alpha)y \|^2 + \alpha(1-\alpha) \|x-y\|^2.
\end{equation}
For the second and third terms on the right hand side of \eqref{equ_thm_57_1}, \vspace{-0.2cm}
\begin{equation}
\begin{array}{ll}\label{equ_thm_57_2}
&\hspace{-0.5cm} \xi \| (  T_1 \varpi_{k}-  T_1 \varpi^*)  - (T_2 T_1 \varpi_{k}-  T_2  T_1 \varpi^*)\|^2_{\Phi} \\
&\qquad+\frac{1-\xi}{\xi}\|  \varpi_{k}-\varpi^*- (T_1 \varpi_{k}- T_1 \varpi^*) \|^2_{\Phi} \\
& =\| (T_1 \varpi_{k}- T_1 \varpi^*)-\xi (T_2 T_1 \varpi_{k}-  T_2  T_1 \varpi^*) \\
&\qquad\qquad  - (1-\xi)(\varpi_{k}-\varpi^*)  \|^2_{\Phi}\\
&\qquad+ \xi(1-\xi)\| (  T_1 \varpi_{k}-  T_1 \varpi^*) - (T_2 T_1 \varpi_{k}-  T_2  T_1 \varpi^*)\\
&\qquad \qquad -(T_1 \varpi_{k}- T_1 \varpi^*)+(\varpi_{k}-\varpi^*) \|_{\Phi}^2\\
&\geq \xi (1-\xi) \| (\varpi_{k}-\varpi^*)-(T_2 T_1 \varpi_{k}-  T_2  T_1 \varpi^*)\|^2_{\Phi}\\
&= \xi(1-\xi) \|\varpi_{k}-T_2 T_1 \varpi_{k} \|^2_{\Phi}\\
&= \xi (1-\xi) \| \varpi_{k}-\varpi_{k+1}\|^2_{\Phi}
\end{array}
\end{equation}
where the first equality follows from \eqref{equation_norm} by setting $\alpha=\xi$, $x= (  T_1 \varpi_{k}-  T_1 \varpi^*) - (T_2 T_1 \varpi_{k}-  T_2  T_1 \varpi^*)$ and $y=(T_1 \varpi_{k}- T_1 \varpi^*)-(\varpi_{k}-\varpi^*)$. Combining \eqref{equ_thm_57_1} and \eqref{equ_thm_57_2} yields $\forall k\geq 0$, \vspace{-0.2cm}
\begin{equation}\label{equ_thm_5}
\|\varpi_{k+1}- \varpi^* \|_{\Phi}^2 \leq   \|\varpi_{k}-  \varpi^*  \|^2_{\Phi} - (1-\xi) \| \varpi_{k}-\varpi_{k+1}\|^2_{\Phi}.
\end{equation}
Using \eqref{equ_thm_5} from $0$ to $k$ and adding all $k+1$ inequalities yields \vspace{-0.25cm}
\begin{equation}
\|\varpi_{k+1}- \varpi^* \|_{\Phi}^2 \leq   \|\varpi_{0}-  \varpi^*  \|^2_{\Phi} - (1-\xi) \sum_{l=0}^{k}\| \varpi_{l}-\varpi_{l+1}\|^2_{\Phi}.
\end{equation}
Taking limit as $k\rightarrow \infty$ we have, \vspace{-0.25cm}
$$(1-\xi) \sum_{k=1}^{\infty} \| \varpi_{k}-\varpi_{k+1}\|^2_{\Phi} \leq \|\varpi_{0}-  \varpi^*  \|^2_{\Phi} $$
Since $1-\xi>0$, it follows that $\sum_{k=1}^{\infty} \| \varpi_{k}-\varpi_{k+1}\|^2_{\Phi}$ converges and  $\lim_{k\rightarrow \infty} \varpi_{k}-\varpi_{k+1}=0$ (since $\Phi \succ0$).\\
Since $\{\| \varpi_{k}-  \varpi^* \| \}$ is bounded   and converges, $\{\varpi_k\}$ is a bounded sequence. Thus, there exists a subsequence$\{\varpi_{n_k}\}$  that converges to some $\overline{\varpi}^*$.
Note that  $T_2\circ T_1$ is continuous and single-valued,   because \eqref{com_fix}
is just an equivalent form of 
Algorithm \ref{dal_1}, and the right hand side of Algorithm \ref{dal_1} is continuous. Since 
${\varpi}_{n_{k}+1}=T_2T_1 \varpi_{n_k}$ and  $T_2T_1$ is continuous, and since $\lim_{n_{k}\rightarrow \infty}\varpi_{n_k}-{\varpi}_{n_{k}+1}=0$, passing to limit point, we have $\overline{\varpi}^* =T_2T_1\overline{\varpi}^* $. Therefore, the limit point $\overline{\varpi}^*$  is a fixed point of $T_2 T_1$, or equivalently, $\overline{\varpi}^* \in zer({\bm{\mathfrak{A}}}+{\bm{\mathfrak{B}}})$. Setting $\varpi^*=\overline{\varpi}^*$ in \eqref{equ_thm_5_2}, it follows that   $\{\| \varpi_{k}-  \overline{\varpi}^* \| \}$ is bounded   and converges. Since there exists a subsequence $\{\varpi_{n_k}\}$ that converges to $\overline{\varpi}^*$, it follows that  $\{\| \varpi_{k}-  \overline{\varpi}^* \| \}$ converges to zero. Therefore, $\lim_{k\rightarrow \infty} \varpi_k \rightarrow \overline{\varpi}^*$ as $k \rightarrow \infty$. Furthermore, by  Lemma \ref{lem_fix_approximation}(ii),  this $\overline{\varpi}^* :=col(\mathbf{x}^*,\bm{z}^*,\bm{\lambda}^*) $ satisfies $col(\mathbf{x}^*,\bm{z}^*,\bm{\lambda}^*) \in zer( {\bm{\mathfrak{A}}} +  {\bm{\mathfrak{B}}})$. Invoking Theorem \ref{thm_zeroHAT_is_correct}(i) concludes the proof. 
\hfill $\Box$
\begin{remark}
Note that,  
since  ${\Phi}^{-1}{\bm{\mathfrak{A}}}$ is maximally monotone, convergence of Algorithm \ref{dal_1} could be proved using a simpler argument, based on \cite[Thm. 25.8]{combettes1}, once  ${\Phi}^{-1}{\bm{\mathfrak{B}}}$ is shown to be cocoercive  (fully, not only \textit{restricted}). Such an argument is used in \cite{PavelCDC2018} (for $c\!=\!1$), under a cocoercivity assumption for $\mathcal{R}^T\mathbf{F}$. Here we do not use any monotonicity assumption on $\mathcal{R}^T\mathbf{F}$, but rather only Lipschitz continuity  of $\mathbf{F}$. 
under which only \textit{restricted} cocercivity of ${\Phi}^{-1}{\bm{\mathfrak{B}}}$ is guaranteed. Thus, we cannot directly apply \cite[Thm. 25.8]{combettes1}. Instead, our elementary proof exploits \textit{restricted} nonexpansiveness properties of $T_1$. 
 {We note that the assumptions on $F$ and $\mathbf{F}$ could be relaxed to hold only locally around $x^*$ and $\mathbf{x}^*$, in which case all results become local. We also note that the class of quadratic games satisfies all assumptions globally.} 
Given a globally known parameter  $\delta$, each agent can {\it independently} choose its local step sizes $\tau_i$, $\nu_i$, and $\sigma_i$ with the rule in \eqref{step_size_choice}, such that ${\Phi}$ in \eqref{metric_matrixBis} is  positive definite. In the case of uniform player step-sizes 
 the step-size bound in \eqref{step_size_choice} simplifies to $\tau \leq (\|\bm{A} \quad \bm{L}_\lambda \| + \delta)^{-1}$.
The bounds in  \eqref{step_size_choice} and Theorem \ref{thm_convergence} recover those in \cite{YiPavelCDC1_2017} (set $c\!=\!0$) and those in \cite{PavelCDC2018} (set $c=1$). 
  \end{remark}

\section{Nash-Cournot game over a network}\label{sec:numerical}


In this section we consider a Nash-Cournot game over a network,  as in \cite{shanbhag1}, \cite{cournotgame}, generalized  by introducing additional market capacities constraints or equivalently globally shared coupling affine constraints. This type of network Cournot game appeared also in the numerical studies of \cite{vanderSchaar} (penalty-based algorithm,  reaching a region near the pure penalized NE), and \cite{YiPavelCDC1_2017} (dual-distributed algorithm), both assuming perfect opponents' decision information.  Other practical decision problems in engineering networks can be described by a Nash-Cournot game over a network, e.g. rate control games in communication networks \cite{shanbhag4},  
demand-response games in smart-grid networks \cite{YeHu_TAC_2017}. 

Consider a set of $N$ firms (players/agents)  involved in the production of a homogeneous commodity  that compete over $m$ markets, $M_1,\cdots,M_m$ (Figure \ref{fig_network_cournot_game}). Firm $i$, $i \in \mathcal{N}$ participates in the competition in $n_i$ markets by deciding to produce and deliver $ x_i \in \mathbf{R}^{n_i}$ amount of products to the markets it connects with.  Its production is limited as $x_i \in \Omega_i \subset \mathbf{R}^{n_i}$.
Firm $i$ has a local matrix $A_i\in \mathbb{R}^{m \times n_i}$  (with elements $1$ or $0$) that specifies which markets it participates in. The $j$-th column of $A_i$, 
$[A_i]_{:j}$ has its $k$-th element as $1$ if and only if player $i$ delivers $[x_i]_j$ amount of production to market $M_k$;  all other elements are $0$. 
Therefore,  $A_1,\cdots,A_N$ can be used to specify a  bipartite graph that represents the connections between firms and markets (see Figure \ref{fig_network_cournot_game}). 
Denote  $n=\sum_{i=1}^N n_i$, $ x=[x_i]_{i  \in \mathcal{N}}\in \mathbb{R}^{n}$,  and   $A=[A_1,\cdots,A_N] \in \mathbb{R}^{m\times n}$. Then $A \, x\in \mathbb{R}^{m}= \sum_{i=1}^N A_i x_i$ is the total product supply to all markets, given the action profile $x$ of all firms. Each market $M_k$ has  a maximal capacity $r_k >0$,  so that $Ax \leq r$, where $r=[r_k]_{k=1,\dots,m} \in \mathbf{R}^m$, should be satisfied;  we consider $b_i = \frac{1}{N} r$. Suppose that $P: \mathbf{R}^m \rightarrow \mathbf{R}^m$ is a price vector function that maps the total supply of each market to the corresponding market's price. Each  firm has a local production cost  $c_i(x_i): \Omega_i \rightarrow \mathbf{R}$. Then the local objective function of company (player) $i$ is $J_i(x_i,x_{-i})= c_i(x_i)-P^T(Ax)A_ix_i$, dependent on the other firms' production profile $x_{-i}$.

\begin{figure}
\begin{center}
\begin{tikzpicture}[->,>=stealth',shorten >=0.3pt,auto,node distance=2.1cm,thick,
  rect node/.style={rectangle, ball color={rgb:red,0;green,20;yellow,0},font=\sffamily,inner sep=1pt,outer sep=0pt,minimum size=12pt},
  wave/.style={decorate,decoration={snake,post length=0.1mm,amplitude=0.5mm,segment length=3mm},thick},
  main node/.style={shape=circle, ball color=green!1,text=black,inner sep=1pt,outer sep=0pt,minimum size=12pt},scale=0.85]


  \foreach \place/\i in {{(-3.1,3.1)/1},{(-3.2,2.2)/2},
  {(-3.3,-0.3)/3},
  {(-3.4,-1.4)/4},
  {(-2.5,3.5)/5},
  {(-2.1,0.5)/6},
  {(-0.6,1.2)/7},
  {(-1.3,-0.6)/8},
  {(-1.4,-1.4)/9},
  {(0.1,3.5)/10},
  {(0.6,2.1)/11},
  {(1.8,0.2)/12},
  {(1.3,-1.3)/13},
  {(2.4,3.4)/14},
  {(1.5,1.5)/15},
  {(2.1,-1.1)/16},
  {(3.2,3.2)/17},
  {(3.3,2.3)/18},
  {(3.4,1.4)/19},
  {(3.5,0.5)/20}}
    \node[main node] (a\i) at \place {};

      \node at (-3.1,3.1){\rm \color{black}{$1$}};
      \node at (-3.2,2.2){\rm \color{black}{$2$}};
      \node at (-3.3,-0.3){\rm \color{black}{$3$}};
      \node at (-3.4,-1.4){\rm \color{black}{$4$}};
      \node at (-2.5,3.5){\rm \color{black}{$5$}};
      \node at (-2.1,0.5){\rm \color{black}{$6$}};
      \node at (-0.6,1.2){\rm \color{black}{$7$}};
      \node at (-1.3,-0.6){\rm \color{black}{$8$}};
      \node at (-1.4,-1.4){\rm \color{black}{$9$}};
      \node at (0.1,3.5){\rm \color{black}{${10}$}};
      \node at (0.6,2.1){\rm \color{black}{${11}$}};
      \node at (1.8,0.2){\rm \color{black}{${12}$}};
      \node at (1.3,-1.3){\rm \color{black}{${13}$}};
      \node at(2.4,3.4) {\rm \color{black}{${14}$}};
      \node at (1.5,1.5){\rm \color{black}{${15}$}};
      \node at (2.1,-1.1){\rm \color{black}{${16}$}};
       \node at (3.2,3.2){\rm \color{black}{${17}$}};
      \node at (3.3,2.3){\rm \color{black}{${18}$}};
      \node at(3.4,1.4) {\rm \color{black}{${19}$}};
      \node at (3.5,0.5){\rm \color{black}{${20}$}};

  \foreach \place/\x in {{(-1.2,2.8)/1},{(-3,1)/2},{(-2.2,-1)/3},
    {(0,0)/4}, {(1,0)/5}, {(2,2)/6}, {(2.5,1)/7}}
  \node[rect node] (b\x) at \place {};

      \node at (-1.2,2.8){\rm \color{red}{$M_1$}};
      \node at (-3,1){\rm \color{red}{$M_2$}};
      \node at (-2.2,-1){\rm \color{red}{$M_3$}};
      \node at (0,0){\rm \color{red}{$M_4$}};
      \node at (1,0){\rm \color{red}{$M_5$}};
      \node at (2,2){\rm \color{red}{$M_6$}};
      \node at (2.5,1){\rm \color{red}{$M_7$}};

         \path[->,blue,thick]               (a1) edge (b1);

         \path[->,blue,thick]               (a2) edge (b1);
         \path[->,blue,thick]               (a2) edge (b2);

         \path[->,blue,thick]               (a3) edge (b2);

         \path[->,blue,thick]               (a4) edge (b3);

         \path[->,blue,thick]               (a5) edge (b1);

         \path[->,blue,thick]               (a6) edge (b1);
         \path[->,blue,thick]               (a6) edge (b2);
         \path[->,blue,thick]               (a6) edge (b3);
         \path[->,blue,thick]               (a6) edge (b4);

         \path[->,blue,thick]               (a7) edge (b4);

         \path[->,blue,thick]               (a8) edge (b3);
         \path[->,blue,thick]               (a8) edge (b4);

         \path[->,blue,thick]               (a9) edge (b3);

         \path[->,blue,thick]               (a10) edge (b1);
         \path[->,blue,thick]               (a10) edge (b4);
         \path[->,blue,thick]               (a10) edge (b6);

         \path[->,blue,thick]               (a11) edge (b4);
         \path[->,blue,thick]               (a11) edge (b5);

         \path[->,blue,thick]               (a12) edge (b5);

         \path[->,blue,thick]               (a13) edge (b5);

         \path[->,blue,thick]               (a14) edge (b6);

         \path[->,blue,thick]               (a15) edge (b5);
         \path[->,blue,thick]               (a15) edge (b6);
         \path[->,blue,thick]               (a15) edge (b7);

         \path[->,blue,thick]               (a16) edge (b5);
         \path[->,blue,thick]               (a16) edge (b7);

         \path[->,blue,thick]               (a15) edge (b7);

         \path[->,blue,thick]               (a17) edge (b6);
         \path[->,blue,thick]               (a17) edge (b7);

         \path[->,blue,thick]               (a18) edge (b7);
         \path[->,blue,thick]               (a19) edge (b7);
         \path[->,blue,thick]               (a20) edge (b7);
\end{tikzpicture}
\end{center}
\caption{Network Nash-Cournot game: An edge from $i$ to $M_k$ on this graph implies that agent/firm $i$ participates in Market $M_k$.}\label{fig_network_cournot_game}
\end{figure}
\begin{figure}
\begin{center}
\begin{tikzpicture}[->,>=stealth',shorten >=0.3pt,auto,node distance=2.1cm,thick,
  rect node/.style={rectangle,ball color=blue!10,font=\sffamily,inner sep=1pt,outer sep=0pt,minimum size=12pt},
  wave/.style={decorate,decoration={snake,post length=0.1mm,amplitude=0.5mm,segment length=3mm},thick},
  main node/.style={shape=circle,ball color=green!1,text=black,inner sep=1pt,outer sep=0pt,minimum size=12pt},scale=0.85]


  \foreach \place/\i in {{(-3,2)/1},{(-2,2)/2},
  {(-1,2)/3},
  {(0,2)/4},
  {(1,2)/5},
  {(2,2)/6},
  {(3,2)/7},
  {(3,1)/8},
  {(3,0)/9},
  {(3,-1)/10},
  {(3,-2)/11},
  {(2,-2)/12},
  {(1,-2)/13},
  {(0,-2)/14},
  {(-1,-2)/15},
  {(-2,-2)/16},
  {(-3,-2)/17},
  {(-3,-1)/18},
  {(-3,-0)/19},
  {(-3,1)/20}}
    \node[main node] (a\i) at \place {};

      \node at (-3,2){\rm \color{black}{$1$}};
      \node at (-2,2){\rm \color{black}{$2$}};
      \node at (-1,2){\rm \color{black}{$3$}};
      \node at (0,2){\rm \color{black}{$4$}};
      \node at (1,2){\rm \color{black}{$5$}};
      \node at (2,2){\rm \color{black}{$6$}};
      \node at (3,2){\rm \color{black}{$7$}};
      \node at (3,1){\rm \color{black}{$8$}};
      \node at (3,0){\rm \color{black}{$9$}};
      \node at (3,-1){\rm \color{black}{${10}$}};
      \node at (3,-2){\rm \color{black}{${11}$}};
      \node at (2,-2){\rm \color{black}{${12}$}};
      \node at (1,-2){\rm \color{black}{${13}$}};
      \node at (0,-2) {\rm \color{black}{${14}$}};
      \node at (-1,-2){\rm \color{black}{${15}$}};
      \node at (-2,-2){\rm \color{black}{${16}$}};
      \node at (-3,-2){\rm \color{black}{${17}$}};
      \node at (-3,-1){\rm \color{black}{${18}$}};
      \node at(-3,0) {\rm \color{black}{${19}$}};
      \node at (-3,1){\rm \color{black}{${20}$}};

        \path[-,blue,thick]               (a1) edge (a2);
            \path[-,blue,thick]               (a2) edge (a3);
         \path[-,blue,thick]               (a3) edge (a4);
         \path[-,blue,thick]               (a4) edge (a5);

         \path[-,blue,thick]               (a5) edge (a6);
         \path[-,blue,thick]               (a6) edge (a7);
         \path[-,blue,thick]               (a7) edge (a8);
         \path[-,blue,thick]               (a8) edge (a9);
           \path[-,blue,thick]               (a9) edge (a10);
         \path[-,blue,thick]               (a10) edge (a11);
             \path[-,blue,thick]             (a11) edge (a12);
            \path[-,blue,thick]               (a12) edge (a13);
            \path[-,blue,thick]               (a13) edge (a14);
            \path[-,blue,thick]               (a14) edge (a15);

            \path[-,blue,thick]               (a15) edge (a16);
            \path[-,blue,thick]               (a16) edge (a17);
            \path[-,blue,thick]               (a17) edge (a18);

            \path[-,blue,thick]               (a18) edge (a19);
            \path[-,blue,thick]               (a19) edge (a20);
             \path[-,blue,thick]               (a1) edge (a20);

              \path[-,blue,thick]               (a2) edge (a15);
             \path[-,blue,thick]               (a6) edge (a13);

\end{tikzpicture}
\end{center}
\caption{Communication graph $\mathcal{G}_c$: Firms $i$ and $j$ are able to exchange their local $\mathbf{x}^i$ and $\mathbf{x}^j$ if there exists an edge between them on this graph. }\label{fig_G_c_graph}
\end{figure}

In the classical,  centralized-information setting each firm is assumed to have instantaneous access to the others' actions $x_{-i}$. This may be impractical in a large network of geographically distributed firms, \cite{cournotgame}. For example, consider that Figure \ref{fig_network_cournot_game} depicts a group of $N=20$ firms located in different continents that participate in $m=7$ markets, with no centralized communication system between them.   Since players  are unable to directly observe the actions of all others, they  engage in local, non-strategic information exchange, to mitigate their lack of global, centralized information. Firms/players  
may communicate with a local subset of neighbouring firms in a peer-to-peer manner,  via some underlying \textit{communication infrastructure}, hence a distributed partial-information setting. A communication network formed $\mathcal{G}_c$ between the firms prescribes how they communicate locally their production decision,  \cite{shanbhag1}.  In this situation, the communication network is formed by the players who are viewed as the nodes in the network.  In this example, we consider that most of the communication is between firms on the same continent, with one or at most two firms in each continent having a direct connection to another firm on another continent. One such instance of the communication network  $\mathcal{G}_c$ is  shown in Figure \ref{fig_G_c_graph}. 
Firms $i$ and $j$ are able to exchange their local variables if there exists an edge between them on this graph. 
Various other topologies can be considered, with different connectivity.  

\vspace{-0.2cm}
\subsection{Assumptions}

Consider that firm $i$'s  production cost  is a strongly convex, quadratic function $c_i(x_i)=x_i^T Q_i x_i+q_i^T x_i$ with
$Q_i\in \mathbf{R}^{n_i \times n_i}$ symmetric and $Q_i \succ 0$ and $q_i\in \mathbf{R}^{n_i}$. Consider that market $M_k$'s price is a linear  function of the total commodity amount supplied to it, $p_k(x)= \bar{P}_k -  \chi_k [Ax]_{k}$ (known as a linear inverse demand function) with $\bar{P}_k, \chi_k \!>\!0$. We show that all assumptions required for Algorithm \ref{dal_1} are satisfied. 

Denote $P=[p_k]_{k=1,m}: \mathbf{R}^n \rightarrow \mathbf{R}^m$,  $\bar{P}=[\bar{P}_k]_{k=1,m} \in \mathbf{R}^m$, 
$\Xi=diag ([\chi_k]_{k=1,m}) \in \mathbf{R}^{m\times m}$, so that $P=\bar{P}-\Xi Ax$ is the vector price function, and $P^TA_ix_i$ is the payoff of firm $i$ obtained by selling  $x_i$ amount  to the markets that it  connects with. Thus, the objective function of player $i$ is,   \vspace{-0.15cm}
\begin{eqnarray}
&J_i(x_i,x_{-i}) =c_i(x_i)- (\bar{P}-\Xi A x)^T A_ix_i, \, \text{ and }   \label{network_cournot_game_function} \\
&\hspace{-0.5cm} \nabla_i J_i(x)= \nabla c_i(x_i) \!+\! A_i^T \Xi A_i x_i \!-\!  A_i^T  (\bar{P}\!-\!\Xi A x), \label{Fi_ex1_0} 
\end{eqnarray}
where $ \nabla c_i(x_i) = 2 Q_i x_i + q_i$. Then, with $x=[x_i]_{i \in \mathcal{N}}$, $A=[A_1, \cdots A_N]$,  $\nabla c(x)  = 2 diag ([Q_i]_{i \in \mathcal{N}}) \, x + [q_i]_{i \in \mathcal{N}}$, we can write $F(x) = [\nabla_iJ_i(x)]_{i\in \mathcal{N}}$  compactly as 
compactly as,  \vspace{-0.15cm}
\begin{equation}\label{F_ex1}
F(x)= Q \, x + h, \quad \text{ where } \, Q := \Sigma +  A^T \Xi A, 
\end{equation}
$\Sigma := diag \big ([2 Q_i \!+\!A_i^T \Xi A_i]_{i \in \mathcal{N}} \big)$ and  
$h := [q_i]_{i \in \mathcal{N}}  \!-\! A^T\bar{P}$. 
Since   $Q_i\succ 0$ and $\Xi \succ 0$ (by  $\chi_k >0$), it follows that   $\Sigma \succ 0$ and $ A^T \Xi A \succeq 0$, hence $Q \succ 0$. Thus,  $F(x)$ is strongly monotone ($\mu = s_{min} (Q) >0$), 
and Lipschitz continuous with   $\theta_0 = \|Q\|$ where $\|Q\| = \sigma_{\max} (Q)$, and Assumptions ~\ref{assum1},  \ref{strgmon_Fassump} are satisfied.
 
{We show next that Assumption~\ref{Lipchitz_boldFassump} is satisfied and $\theta=\theta_0$. \\
Using \ \eqref{Fi_ex1_0} with $x$ replaced by $\mathbf{x}^i = (x_i,\mathbf{x}^i_{-i})$, yields  \vspace{-0.15cm}
$$
\nabla_i J_i(\mathbf{x}^i) = \nabla c_i(x_i) + A_i^T \Xi A_i x_i-  A_i^T \, (\bar{P}-\Xi A \mathbf{x}^i).
$$
Thus,  $\mathbf{F}(\mathbf{x}):=[\nabla_iJ_i(\mathbf{x}^i)]_{i\in \mathcal{N}}$, \eqref{pseudogradient_EKT}, is given by 
$\mathbf{F}(\mathbf{x}) \!= \!\nabla c (x) \! + \! diag \big ([A_i^T \Xi A_i]_{i\in \mathcal{N}} \big) x \!  - \!  A^T\bar{P}
\! + \! diag \big( [A_i^T \Xi A]_{i \in \mathcal{N}}\big) \mathbf{x} 
$
where  $\mathbf{x} = [\mathbf{x}^i]_{i\in \mathcal{N}}$. Combining the first three terms (with $\nabla c (x) $,  $x=\mathcal{R} \mathbf{x}$, $\Sigma$ and $h$ defined above), and  $diag \big ( [\mathcal{R}_i A^T \Xi A]_{i \in \mathcal{N}} \! \big) =\mathcal{R}( I_N \otimes  A^T \Xi A)=$ for the last term, we can write 
$\mathbf{F}(\mathbf{x}) =  \Sigma  \mathcal{R} \mathbf{x} + h+  \mathcal{R} (I_N \otimes  A^T \Xi A) \, \mathbf{x}.
$ 
Since $\Sigma$ is block-diagonal,  $\Sigma  \mathcal{R}= \mathcal{R} \! (I_N \otimes \Sigma)$, so that 
$\mathbf{F}(\mathbf{x}) \! = \!   \mathcal{R}  \big (I_N \otimes  (\Sigma  \!+ \! A^T \Xi A)  \big ) \, \mathbf{x} + h $, which, with $Q$ as in \eqref{F_ex1}, is \vspace{-0.15cm} 
 \begin{equation}\label{boldF_ex1_final}
 \mathbf{F}(\mathbf{x})=  \mathcal{R}  (I_N \otimes Q)\,\mathbf{x} + h.  
\end{equation}  
Hence, $\|\mathbf{F}(\mathbf{x}) -\mathbf{F}(\mathbf{y}) \| \leq \| \mathcal{R} (I_N \otimes Q)\|  \| \mathbf{x} - \mathbf{y} \| \leq \| \mathcal{R}\| \| I_N \otimes Q\|   \| \mathbf{x} - \mathbf{y} \| \leq  \|Q\|    \| \mathbf{x} - \mathbf{y} \| $, since $ \| \mathcal{R}\| =1$ 
and $\| I_N \otimes Q\| = \| Q \|$. Hence, $\mathbf{F}(\mathbf{x})$ is Lipschitz and Assumption \ref{Lipchitz_boldFassump} holds with $\theta =\theta_0= \|Q\|$. 
}

\vspace{-0.2cm}
\subsection{Numerical results}

\begin{figure}
  \centering
  \includegraphics[width=3.5in]{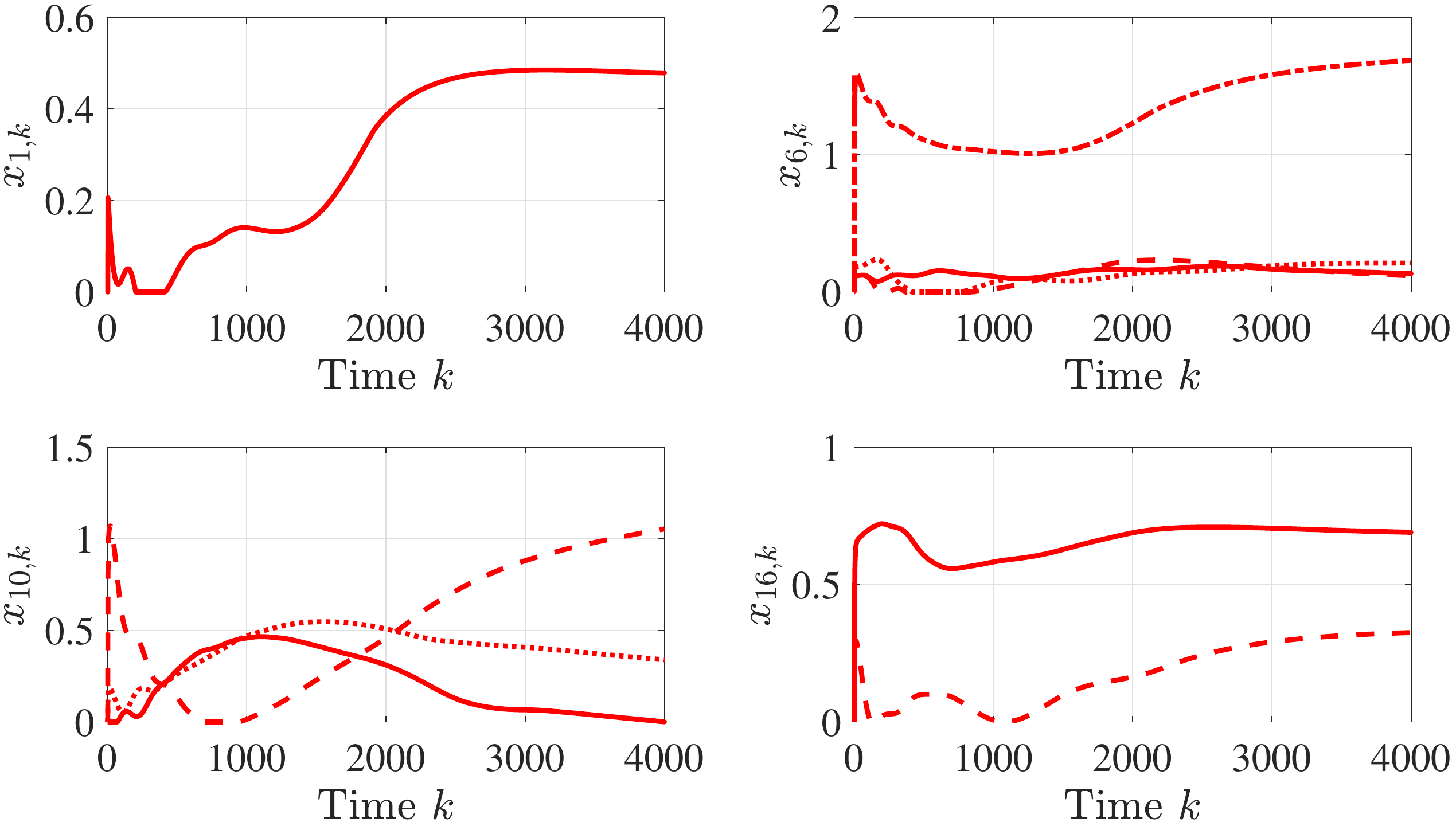}\\
  \caption{The trajectories of local decisions $x_{i,k}$ for selected players ($c=100$) }\label{fig_sim_3_100}  
\end{figure}

\begin{figure}
  \centering
  \includegraphics[width=3.5in]{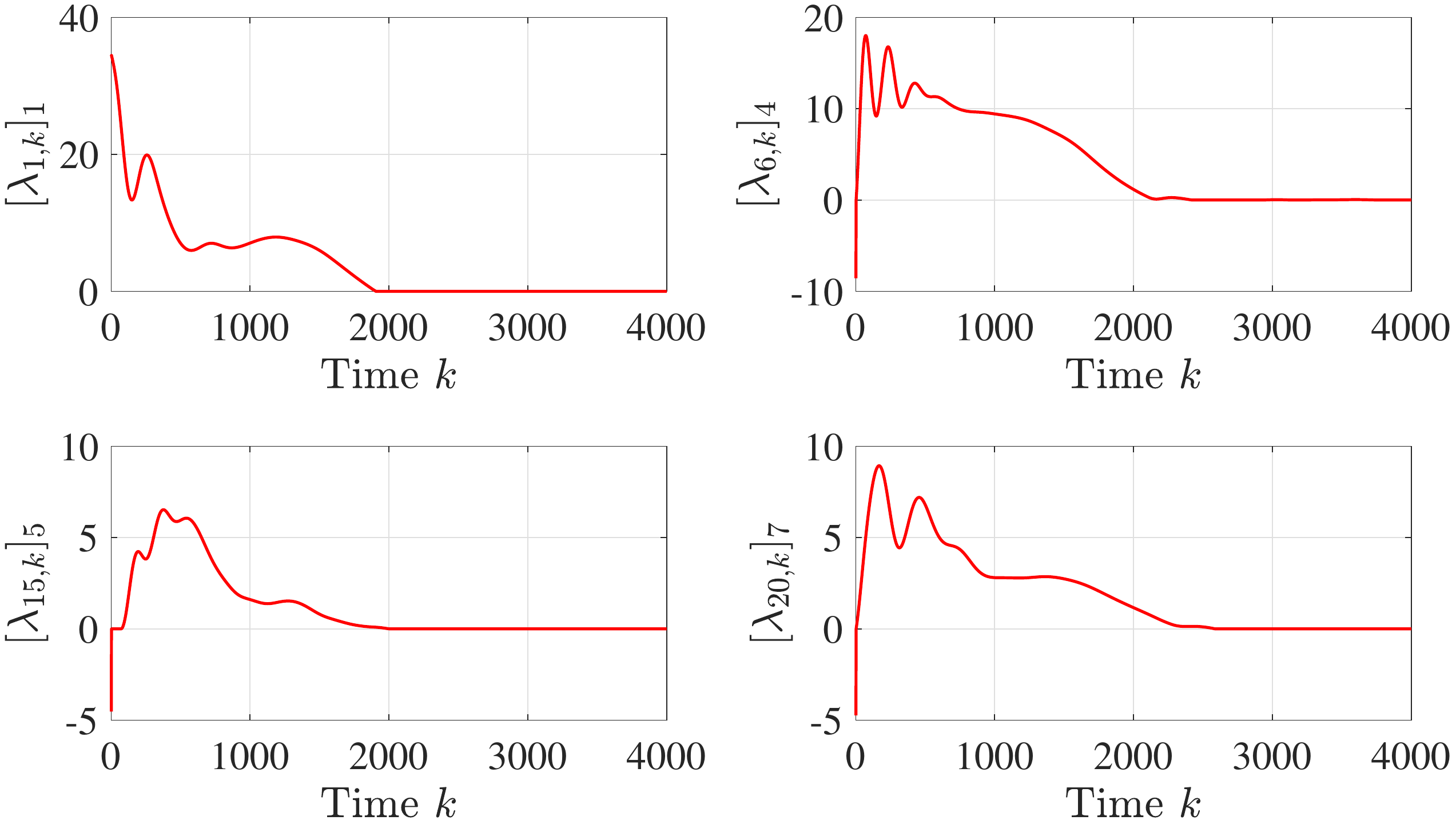}\\
  \caption{The trajectories of dual variables $\lambda_{i,k}$  for selected players ($c=100$)  }\label{fig_sim_2_100}  
\end{figure}

Consider Figure \ref{fig_network_cournot_game} ($N\!=\!20$ and $m\!=\!7$). For the communication graph $\mathcal{G}_c$  in Fig. \ref{fig_G_c_graph}, with weighted adjacency matrix $W=[w_{ij}]$ having all its nonzero elements as $1$, $\lambda_2(L) = 0.102$. We consider $ \bm{0}\leq x_i \leq X_i $ where  each component of $X_i$ is randomly drawn from $(5,10)$. Each market $M_k$ has a maximal capacity  of $r_k$, where  $r_k$ is  randomly  drawn from $(1,2)$. Player $i$'s objective function is taken as \eqref{network_cournot_game_function}, where in  $c_i(x_i)=x_i^T Q_i x_i+q_i^T x_i$, $Q_i$ is diagonal with its entries   
randomly drawn from $(1,8)$, and each component of $q_i$ is randomly from $(1,2)$. In the price function $ P=\bar{P}-\Xi Ax$, $\bar{P}_k$ and $\chi_k$ are  randomly drawn from $(10,20)$ and $(1,3)$, respectively. For  $Q$, \eqref{F_ex1}, this yields  $s_{\min}(Q) = 1.001=\mu$,  $s_{\max}(Q) = 2.09=\theta$, and 
  the lower bound for $c$ in Lemma \ref{boldRF_plus_cL_strgmon_lemma} is $62.5$. We set   $c=100$ and $\tau_i=0.003$, $\nu_i=0.02$ and $\sigma_i=0.003$ for all firms, which  satisfy \eqref{step_size_choice}. 
The results for the implementation of Algorithm \ref{dal_1}  in the partial-decision setting, where the primal  and dual variables 
are exchanged over  the sparsely connected graph $\mathcal{G}_c$ (Fig.~\ref{fig_G_c_graph}),  are shown in Fig. \ref{fig_sim_3_100}, \ref{fig_sim_2_100}. 

\begin{figure}
  \centering
  \includegraphics[width=3.5in]{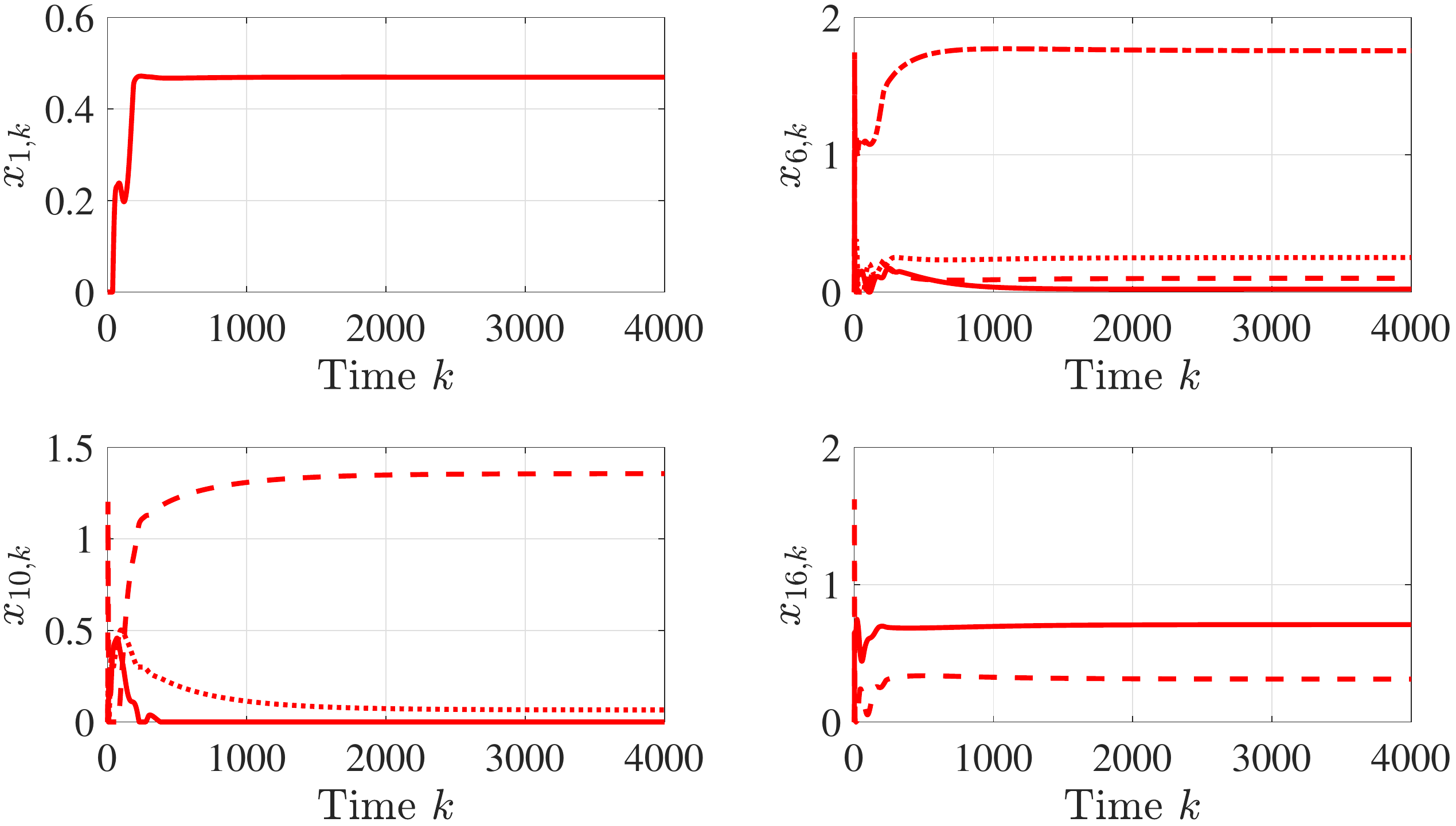}\\
  \caption{The trajectories of local decisions $x_{i,k}$ for selected players ($c=10$) }\label{fig_sim_3_10}   
\end{figure}

The algorithm converges to the same GNE found by \cite{YiPavelCDC1_2017}, with a comparable rate of convergence. Unlike \cite{YiPavelCDC1_2017},  here  the primal variables are not perfectly known, but estimated. 
In fact, the lower bound for $c$ in Lemma \ref{boldRF_plus_cL_strgmon_lemma} is quite conservative. For comparison, we set $c=10$ and  increase the step-sizes ten times (still satisfying \eqref{step_size_choice}). The simulation results for the same initial conditions are shown in  Fig. \ref{fig_sim_3_10} and indicate fast convergence. 


\vspace{-0.3cm}
\section{Conclusions}\label{sec_concluding}

In this paper, we  considered a partial-decision information setting and  we proposed a fully-distributed, primal-dual algorithm  for computation of a variational GNE in noncooperative games with globally-shared affine coupling constraints. The algorithm is motivated  by a forward-backward splitting method for finding zeros of a sum of doubly-augmented monotone operators. We proved its convergence 
 with  fixed step-sizes over any connected graph,  
 by  leveraging monotone operator-splitting techniques. In future work we will consider mechanism designs so that players faithfully report their variables. 



\vspace{0.2cm}
\noindent {\bf Acknowledgements:} The author would like to acknowledge Bolin Gao's help in performing the numerical experiments. 

\vspace{-0.3cm}

\appendix

{\bf Proof  of Lemma \ref{lem_ALG_compact}:}
With 
$x_{i,k}\!=\!\mathbf{x}^i_{i,k}\!=\!\mathcal{R}_i  \mathbf{x}^i_k$,  $\mathbf{x}^j_{i,k} \!=\!\mathcal{R}_i  \mathbf{x}^j_k$,  the first line of Algorithm \ref{dal_1} is written  as 
$
x_{i,k+1} \!    = \!  P_{\Omega_i}\big[\!x_{i,k} \!-\!\tau_i (\! \nabla_{x_i} J_i(x_{i,k},\mathbf{x}^i_{-i,k}) \!+\!A_i^T \lambda_{i,k}   
\! + \! c \mathcal{R}_i \sum_{j\in \mathcal{N}_i} w_{ij} (\mathbf{x}^i_{k} \!-\! \mathbf{x}^j_{k}) )\!\big].
$ 
Then,  for all players  in stacked notation, with  $\mathbf{F}(\mathbf{x})$, \eqref{pseudogradient_EKT}, $\bm{A}$, $\mathbf{L}_{x} $, $\mathcal{R}$ and   $P_{\prod_{i=1}^N \Omega_i}(x) = col(P_{\Omega_i}(x_i))_{i \in \mathcal{N}}$,  \cite[Prop. 23.16]{combettes1},   the $x$-update can be written as \eqref{Alg1BLOCK_1}.  
With $\mathbf{x}^j_{-i,k} \!=\!\mathcal{S}_i \, \mathbf{x}^j_k$, the second line in  Algorithm \ref{dal_1} is  
$
\mathcal{S}_i \, \mathbf{x}^i_{k+1}  \!=\! \mathcal{S}_i \, \mathbf{x}^i_k  - \tau_i {c} \mathcal{S}_i  \sum_{j\in \mathcal{N}_i} w_{ij}  (\mathbf{x}^i_{k}- \mathbf{x}^j_{k}).
$   
In stacked notation, with  $ col(\mathbf{x}^i_{-i})_{i \in \mathcal{N}}= \mathcal{S}\mathbf{x}$, this  gives  \eqref{Alg1BLOCK_2}. Proceeding similarly for the 3$^{rd}$ and 4$^{th}$  line of Algorithm \ref{dal_1} yields \eqref{Alg1BLOCK_3} and \eqref{Alg1BLOCK_4}. 
\hfill $\Box$

{\bf Proof of Lemma \ref{lem_fix_approximation}:} 
(i) We show that  expanding \eqref{compact_operator_1} for  $\varpi_k\!=\!col({\mathbf{x}}_k, {\bm{z}}_k,{\bm{\lambda}}_k)$, we obtain  \eqref{Alg1BLOCK_1}-\eqref{Alg1BLOCK_4} in  Lemma \ref{lem_ALG_compact}, hence Algorithm \ref{dal_1}. First, recall $\bm{\tau}_x\! =\! diag((\tau_i I_{n_i})_{i\in \mathcal{N}})$.  
By using  $\mathcal{R} \!=\! diag((\mathcal{R}_i)_{i\in \mathcal{N}}) $ and \eqref{eq:actualStratREMatrix},   it can be shown that 
 $ \bm{\tau}_x \mathcal{R}\!=\!\mathcal{R}\bm{\tau} $,  where  $\bm{\tau} \!=\! diag((\tau_i I_n)_{i\in \mathcal{N}})$, \eqref{metric_matrixBis}. Similarly, $\bm{\tau}_s \mathcal{S}\!=\!\mathcal{S}\bm{\tau} $. Hence, with $\mathcal{R}\mathcal{R}^T \!=\!  I_{n}$ and $\mathcal{R}^T\mathcal{R}\!+\! \mathcal{S}^T\mathcal{S}\!=\! I_{Nn}$ (cf.  \eqref{eq:IdRSMatrix}), this yields  \vspace{-0.2cm}
 $$\bm{\tau}_x\!=\!\mathcal{R}\bm{\tau}\mathcal{R}^T \quad \text{and} \quad  
 \mathcal{R}^T\bm{\tau}_x\mathcal{R}\! +\!\mathcal{S}^T\bm{\tau}_s\mathcal{S}\!=\! \bm{\tau}.$$ 

Using ${\bm{\mathfrak{A}}}$, ${\bm{\mathfrak{B}}}$,  \eqref{op_hat_A}, $\Phi$, \eqref{metric_matrixBis},  the  update for ${\mathbf{x}}_k$  in \eqref{compact_operator_1} is    \vspace{-0.2cm}
\begin{equation}\label{eq_c_11}
\begin{array}{ll}
 -  \mathcal{R}^T  \mathbf{F}(\mathbf{x}_{k})  & - {c}\mathbf{L}_x \,  \mathbf{x}_k \in  \mathcal{R}^T N_{\Omega}( \mathcal{R}\mathbf{x}_{k+1})+ \mathcal{R}^T \bm{A}^T \bm{\lambda}_{k+1} \\
&  + \bm{\tau}^{-1}(\mathbf{x}_{k+1}- \mathbf{x}_{k})     - \mathcal{R}^T \bm{A}^T (\bm{\lambda}_{k+1}- \bm{\lambda}_{k})
\end{array}
\end{equation}
which  means that for some $v \in N_{\Omega}( \mathcal{R}\mathbf{x}_{k+1})=N_{\Omega}( x_{k+1})$, \vspace{-0.2cm} 
\begin{equation}\label{eq_c_11_NEW}
 -  \bm{\tau} \mathcal{R}^T  \mathbf{F}(\mathbf{x}_{k}) \! - \! \bm{\tau} {c}\mathbf{L}_x \,  \mathbf{x}_k \!= \!  \bm{\tau} \mathcal{R}^T ( v \!+\!   \bm{A}^T \bm{\lambda}_{k} )
  \!+\! \mathbf{x}_{k+1} \!-\! \mathbf{x}_{k}.
\end{equation}
Premultiplying \eqref{eq_c_11_NEW} by $\mathcal{R}$ 
 yields \vspace{-0.2cm}
$$- \!\mathcal{R} \bm{\tau} \!\mathcal{R}^T  \mathbf{F}(\mathbf{x}_{k})\!  -\! \mathcal{R} \bm{\tau} {c}\mathbf{L}_x   \mathbf{x}_k 
\!=\!  \mathcal{R} \bm{\tau}\! \mathcal{R}^T ( v +   \bm{A}^T \bm{\lambda}_{k} )
 \!+\! \mathcal{R} (\mathbf{x}_{k+1}- \mathbf{x}_{k})
$$ 
With  $\mathcal{R}\bm{\tau}\mathcal{R}^T \!= \!\bm{\tau}_x$, $\mathcal{R}\bm{\tau} \!= \!\bm{\tau}_x \mathcal{R}$,  and  $  \mathcal{R} \mathbf{x} \!=x$,   this yields\vspace{-0.2cm}
$$- \bm{\tau}_x   \mathbf{F}(\mathbf{x}_{k})  - \bm{\tau}_x {c}\mathcal{R}  \mathbf{L}_x \,  \mathbf{x}_k 
=  \bm{\tau}_x ( v +   \bm{A}^T \bm{\lambda}_{k} )  +  x_{k+1}- x_{k},
$$ for $v \in N_{\Omega}( x_{k+1})$, hence \vspace{-0.2cm} %
$$
x_{k}-\bm{\tau}_x (  \mathbf{F}(\mathbf{x}_{k}) +\bm{A}^T \bm{\lambda}_{k} +{c} \mathcal{R} \mathbf{L}_x \,  \mathbf{x}_k)  \in  x_{k+1}+ \bm{\tau}_x N_{\Omega}(x_{k+1}).
$$
Note that   $\bm{\tau}_x N_{\Omega}(x)= \prod_{i=1}^N \tau_i N_{\Omega_i}(x_i) =N_{\Omega}(x)$, since $\tau_i >0$, $N_{\Omega}$ is a cone and $N_{\Omega}(x) = \prod_{i=1}^N N_{\Omega_i}(x_i)$. With  $P_{\Omega}(x)=({\rm Id}+N_{\Omega})^{-1}(x)$, the foregoing gives  \eqref{Alg1BLOCK_1}.  

Premultiplying \eqref{eq_c_11_NEW} by $\mathcal{S}$, using $\mathcal{S}\bm{\tau} \!\mathcal{R}^T\!=\!\bm{\tau}_s \mathcal{S} \!\mathcal{R}^T $, $\mathcal{S}\mathcal{R}^T\!=\!\bm{0}$, \eqref{eq:IdRSMatrix} yields 
$ -  \mathcal{S} \bm{\tau}{c}\mathbf{L}_x \, \mathbf{x}_k \!=\!  \mathcal{S} (\mathbf{x}_{k+1}   \!-\! \mathbf{x}_{k})$. With  $\mathcal{S}\bm{\tau}\! =\! \bm{\tau}_s \mathcal{S}$, this is  \eqref{Alg1BLOCK_2}. The reverse direction, from \eqref{Alg1BLOCK_1},\eqref{Alg1BLOCK_2} to \eqref{eq_c_11} can be shown similarly, with $\mathcal{R}^T\mathcal{R}\!+\! \mathcal{S}^T\mathcal{S}\!=\! I_{Nn}$, $\!\mathcal{R}^T\bm{\tau}_x\mathcal{R}\! +\!\mathcal{S}^T\bm{\tau}_s\mathcal{S}\!=\! \bm{\tau}$. Thus the   ${\mathbf{x}}_k$ update  in \eqref{compact_operator_1} is equivalent to  \eqref{Alg1BLOCK_1}, \eqref{Alg1BLOCK_2}.

The $\bm{z}_k$ update in \eqref{compact_operator_1} is 
$\! \bm{0}
 =            			\!-\mathbf{L} _\lambda \bm{\lambda}_{k+1}
                  +   \! \bm{\nu}^{-1} (\bm{z}_{k+1} -\bm{z}_{k}) +  \!\mathbf{L}_\lambda (\bm{\lambda}_{k+1} -  \bm{\lambda}_{k}),
$,  
which is   \eqref{Alg1BLOCK_3}. The $\bm{\lambda}_k$ update in \eqref{compact_operator_1} is  \vspace{-0.2cm} 
\begin{align*}
 -[\mathbf{L}_{\lambda} \, \bm{\lambda}_{k} +\bm{b}]  \in N_{\bm{R}^{Nm}_{+}}(\bm{\lambda}_{k+1})-\bm{A} \mathcal{R} \mathbf{x}_{k+1}+\mathbf{L}_{\lambda} \,\bm{z}_{k+1} \qquad \qquad \nonumber \\
-\bm{A} \mathcal{R}  (\mathbf{x}_{k+1}-\mathbf{x}_{k})+\mathbf{L}_{\lambda} \,(\bm{z}_{k+1}-\bm{z}_k)+\bm{\sigma}^{-1}(\bm{\lambda}_{k+1}-\bm{\lambda}_{k}).\,\,\,\, \nonumber
 \end{align*} 
With    $x =  \mathcal{R} \mathbf{x}$ this reduces to \vspace{-0.2cm}
\begin{align*}
 \bm{\lambda}_{k} +\bm{\sigma}[    \bm{A} (2x_{k+1} -x_{k}) - \bm{b} 
 -\mathbf{L}_{\lambda} \,(2\bm{z}_{k+1}-\bm{z}_k) -\mathbf{L}_{\lambda} \, \bm{\lambda}_{k} ] \nonumber \\
 \in N_{\bm{R}_{+}^{Nm}}(\bm{\lambda}_{k+1})+ \bm{\lambda}_{k+1},
 \end{align*} 
which is   \eqref{Alg1BLOCK_4}, by $\bm{\sigma}^{-1}N_{\bm{R}^{Nm}_{+}}(\bm{\lambda})$ $\!=\!N_{\bm{R}^{Nm}_{+}}(\bm{\lambda})$ and  
$P_{\bm{R}_{+}^{Nm}}(\bm{\lambda})\!=\!({\rm Id}\!+\!N_{\bm{R}_{+}^{Nm}})^{-1}$. Thus,   \eqref{Alg1BLOCK_1}-\eqref{Alg1BLOCK_4} are equivalent to \eqref{compact_operator_1}. 

Since $\Phi \succ 0$,   \eqref{compact_operator_1} is equivalent to 
\vspace{-0.2cm}
\begin{equation}\label{compact_operator_2}
({\rm Id}-{\Phi}^{-1}{\bm{\mathfrak{B}}})(\varpi_{k}) \in ({\rm Id}+{\Phi}^{-1}{\bm{\mathfrak{A}}})(\varpi_{k+1}).
\end{equation}
In turn,  since ${\Phi}^{-1} {\bm{\mathfrak{A}}}$ is maximally monotone,   it follows that   $({\rm Id}+{\Phi}^{-1} {\bm{\mathfrak{A}}})^{-1}$ is single-valued (by \cite[Prop.  23.7]{combettes1}), so that  \eqref{compact_operator_2} is equivalently written as  
 \eqref{com_fix}.  

(ii)  Suppose that Algorithm \ref{dal_1}, or \eqref{com_fix},  has a limit point $\overline{\varpi} = col(\overline{\mathbf{x}}, \overline{\bm{z}},\overline{\bm{\lambda}})$. By the continuity of the right-hand-side of Algorithm \ref{dal_1},  the following  equivalences hold, \vspace{-0.2cm} 
\begin{equation}\label{eq:fixpt_T2T1}
\overline{\varpi}=({\rm Id}+{\Phi}^{-1} {\bm{\mathfrak{A}}})^{-1}\circ ({\rm Id}- {\Phi}^{-1} {\bm{\mathfrak{B}}}) \overline{\varpi}{ : = T_2 \circ T_1  \overline{\varpi}}
\end{equation}
$$ \Leftrightarrow ({\rm Id}-{\Phi}^{-1}{\bm{\mathfrak{B}}})(\overline{\varpi} ) \in ({\rm Id}+{\Phi}^{-1}{\bm{\mathfrak{A}}})(\overline{\varpi} )$$
$$ \Leftrightarrow \bm{0}  \in ({\Phi}^{-1}{\bm{\mathfrak{A}}}+{\Phi}^{-1}{\bm{\mathfrak{B}}})(\overline{\varpi} )$$
$$
 \Leftrightarrow 
 \bm{0}  \in ({\bm{\mathfrak{A}}}+{\bm{\mathfrak{B}}})(\overline{\varpi} ).
$$ Hence, $\overline{\varpi} \in zer ({\bm{\mathfrak{A}}}+{\bm{\mathfrak{B}}})$ and  is a fixed point of $T_2 \circ T_1$.
\hfill $\Box$

{{\bf Proof of Lemma \ref{boldRF_plus_cL_strgmon_lemma}}}:  
Recall the decomposition of  $\mathbf{R}^{Nn}$ 
into $\mathbf{R}^{Nn} \! \!=\! \mathbf{E}_x \oplus \mathbf{E}_x^\perp$ where $\mathbf{E}_x \! =\! Null(\mathbf{L}_x)$ is the  estimate consensus subspace and $\mathbf{E}_x^\perp \!=\! Null(\mathbf{L}_x)^{\perp}$ its orthogonal complement. 
Any $\mathbf{x}  \in  \mathbf{R}^{Nn}$ can be decomposed as $\mathbf{x}= \mathbf{x}^{\|} + \mathbf{x}^{\perp}$, with $ \mathbf{x}^{\|} \in  \mathbf{E}_x$, $\mathbf{x}^{\perp} \in \mathbf{E}_x ^{\perp}$ and  $ ( \mathbf{x}^{\|})^T\mathbf{x}^{\perp} =0$, 
by using two projection matrices   $P_{\|} = \frac{1}{N}  \mathbf{1}_N \otimes  \mathbf{1}_N^T \otimes I_n$, $P_{\perp} = I_{Nn} -  \frac{1}{N}  \mathbf{1}_N \otimes  \mathbf{1}_N^T \otimes I_{n}$.  
Thus $\mathbf{x}^{\|} = \mathbf{1}_N \otimes x$, for some  $x \in  \mathbf{R}^{n}$, so that  $\mathbf{L}_x \mathbf{x}^{\|} =0$, and $min_{\mathbf{x}^{\perp}\in\mathbf{E}_x ^{\perp}} (\mathbf{x}^{\perp})^T \mathbf{L}_x \mathbf{x}^{\perp} = s_2(L) \| \mathbf{x}^{\perp} \|^2$, where $ s_2(L) >0$.
Any $\mathbf{x'} \in \mathbf{E}_x$ 
can be written as $\mathbf{x'}= \mathbf{1}_N \otimes x'$, 
for some  $x' \in  \mathbf{R}^{n}$, such that  $\mathbf{L}_x \mathbf{x'}=0$. 
Then, using  $\mathbf{x} = \mathbf{x}^{\|} + \mathbf{x}^{\perp}$, $\mathbf{F}(\mathbf{x}^{\|}) = F(x)$, $\mathbf{F}(\mathbf{x'}) = F(x')$, $\mathcal{R} \mathbf{x}^{\|} = x$, $\mathcal{R} \mathbf{x'} = x'$, yields  \vspace{-0.15cm}
\begin{eqnarray}\label{eqboldRF_plus_cL_strgmon_1}
&&\hspace{-0.3cm}(\!\mathbf{x}\!-\! \mathbf{x'})^T \! \big ( \mathcal{R}^T (\mathbf{F}(\mathbf{x})\!-\!\mathbf{F}(\mathbf{x'}))\!+\!c \mathbf{L}_x(\mathbf{x} \!-\! \mathbf{x'}) \big )\!  \nonumber \\
&&= (\!\mathbf{x}^{\|}\!-\! \mathbf{x'})^T\mathcal{R}^T \big (\mathbf{F}(\mathbf{x})\!-\!\mathbf{F}(\mathbf{x}^{\|})\! +\! \mathbf{F}(\mathbf{x}^{\|})\!  -\!\mathbf{F}(\mathbf{x'}) \big)\! \nonumber \\
&&\hspace{0.2cm} + (\mathbf{x}^{\perp})^T\mathcal{R}^T \big ( \mathbf{F}(\mathbf{x})\!-\!\mathbf{F}(\mathbf{x}^{\|})\! +\! \mathbf{F}(\mathbf{x}^{\|})\!  -\!\mathbf{F}(\mathbf{x'}) \big )\! \nonumber \\
&&\hspace{0.2cm} + c  (\mathbf{x}^{\|}\!+\!\mathbf{x}^{\perp}-\! \mathbf{x'})^T\mathbf{L}_x(\mathbf{x}^{\|} + \mathbf{x}^{\perp} \!-\! \mathbf{x'}) \big ) \! \nonumber \\
&&= (x-x')^T \!\big(\mathbf{F}(\mathbf{x})\!-\!\mathbf{F}(\mathbf{x}^{\|}) \big )\! +  (x-x')^T \! \big ( F(x)\!  -F(x') \big )\! \nonumber \\
&&\hspace{0.2cm} + (\!\mathbf{x}^{\perp})^T \mathcal{R}^T \big (\mathbf{F}(\mathbf{x})\!-\!\mathbf{F}(\mathbf{x}^{\|}) \big )\! + \! (\!\mathbf{x}^{\perp})^T \mathcal{R}^T  \big (F(x)\!  -\! F(x') \big )\! \nonumber \\
&&\hspace{0.2cm} +c (\mathbf{x}^{\perp})^T \mathbf{L}_x \mathbf{x}^{\perp} 
\end{eqnarray}
Using strong monotonicity of $F$ (by Assumption \ref{strgmon_Fassump}) for the second term, and properties of $ \mathbf{L}_x $ on $\mathbf{E}_x ^\perp=Null(\mathbf{L}_x)^{\perp} $ for the fifth term  on the right-hand side we can write,  \vspace{-0.15cm}
\begin{eqnarray}\label{eqboldRF_plus_cL_strgmon_2}
&&\hspace{-0.3cm}(\!\mathbf{x}\!-\! \mathbf{x'})^T \! \big ( \mathcal{R}^T (\mathbf{F}(\mathbf{x})\!-\!\mathbf{F}(\mathbf{x'}))\!+\!c \mathbf{L}_x(\mathbf{x} \!-\! \mathbf{x'}) \big )\!  \nonumber \\
&& \geq (x-x')^T \!\big(\mathbf{F}(\mathbf{x})\!-\!\mathbf{F}(\mathbf{x}^{\|}) \big )\! + \mu \| x-x' \|^2 \nonumber \\
&&\hspace{0.2cm} + (\!\mathbf{x}^{\perp})^T \mathcal{R}^T \big (\mathbf{F}(\mathbf{x})\!-\!\mathbf{F}(\mathbf{x}^{\|}) \big )\! + \! (\!\mathbf{x}^{\perp})^T \mathcal{R}^T  \big (F(x)\!  -\! F(x') \big )\! \nonumber \\
&&\hspace{0.2cm} + c  s_2(L)\| \mathbf{x}^{\perp} \|^2
\end{eqnarray}
We deal with the cross-terms by using $a^T b \geq - \|a\| \|b \|$, for any $a$, $b$, and $\| \mathcal{R}^T b \| \leq \| \mathcal{R}^T\| \|b\| $, so that  
\vspace{-0.25cm}
\begin{eqnarray}\label{eqboldRF_plus_cL_strgmon_3}
&&\hspace{-0.5cm}(\!\mathbf{x}\!-\! \mathbf{x'})^T \! \big ( \mathcal{R}^T (\mathbf{F}(\mathbf{x})\!-\!\mathbf{F}(\mathbf{x'}))\!+\!c \mathbf{L}_x(\mathbf{x} \!-\! \mathbf{x'}) \big )\! \geq	 \nonumber \\
&&\hspace{-0.3cm} - \|x-x'\| \, \| \mathbf{F}(\mathbf{x})\!-\!\mathbf{F}(\mathbf{x}^{\|}) \| \! + \mu \| x-x' \|^2 + c  s_2(L)\| \mathbf{x}^{\perp} \|^2 \nonumber \\
&&\hspace{-0.3cm}- \|\!\mathbf{x}^{\perp} \| \|\mathcal{R}^T\| \| \mathbf{F}(\mathbf{x})\!-\!\mathbf{F}(\mathbf{x}^{\|}) \| \! - \! \| \!\mathbf{x}^{\perp}\| \|\mathcal{R}^T\| 
\| F(x)\!  -\! F(x') \|\! \nonumber 
\end{eqnarray}
Using  $\| \mathcal{R}^T\| = 1$ and  Lipschitz properties of $F$ and $\mathbf{F}$,  \vspace{-0.15cm} 
\begin{eqnarray}\label{eqboldRF_plus_cL_strgmon_4}
&&\hspace{-0.5cm}(\!\mathbf{x}\!-\! \mathbf{x'})^T \! \big ( \mathcal{R}^T (\mathbf{F}(\mathbf{x})\!-\!\mathbf{F}(\mathbf{x'}))\!+\!c \mathbf{L}_x(\mathbf{x} \!-\! \mathbf{x'}) \big )\!  \geq \nonumber \\
&& \hspace{-0.3cm} - \theta \|x-x'\| \, \| \mathbf{x}^{\perp} \| \! + \mu \| x-x' \|^2 - \theta \| \mathbf{x}^{\perp} \|^2 \!  - \theta_0 \! \| \mathbf{x}^{\perp}\| 
\| x\!  -\! y \|\!  \nonumber \\
&&\hspace{-0.3cm}
 + c  s_2(L)\| \mathbf{x}^{\perp} \|^2
\end{eqnarray}
Using $ \|x-x'\| = \frac{1}{\sqrt{N}}  \| \mathbf{x}^{\|} \!-\! \mathbf{x'}\|$ and $\Psi $, \eqref{eqboldRF_plus_cL_strgmon_5}, we can write 
\vspace{-0.25cm}
\begin{eqnarray*}
&&\hspace{-1.8cm}(\!\mathbf{x}\!-\! \mathbf{x'})^T \!\! \big ( \!\mathcal{R}^T (\mathbf{F}(\mathbf{x}\!)\!-\!\mathbf{F}(\mathbf{x'})\!)\!+\!c \mathbf{L}_x(\mathbf{x} \!-\! \mathbf{x'}\!)\! \big )\!  \\
&&\qquad \qquad  \geq \!
\!\left [ \! \begin{array}{c}
\! \| \mathbf{x} ^{\|}\!-\! \mathbf{x'} \| \!\\
  \| \mathbf{x}^{\perp} \|
 \end{array} 
\!  \right ]^{\!T}
\! \!\! \!\! \Psi  \! \left [ 
\! \begin{array}{c}
 \! \| \mathbf{x} ^{\|}\!-\! \mathbf{x'} \|\! \\
  \| \mathbf{x}^{\perp} \|
 \end{array} 
 \! \right ]
\end{eqnarray*}
Thus,  for any $\mathbf{x}$ and  any $\mathbf{x'} \in \mathbf{E}_x$, 
$
(\!\mathbf{x}\!-\! \mathbf{x'})^T \! \big ( \mathcal{R}^T \mathbf{F}(\mathbf{x}\!)\!-\! \mathcal{R}^T \mathbf{F}(\mathbf{x'})\!+\!c \mathbf{L}_x(\mathbf{x} \!-\! \mathbf{x'}\!) \big )\!  \geq \! \bar{\mu} 
\left (
 \| \mathbf{x} ^{\|}\!-\! \mathbf{x'} \|^2 + 
  \| \mathbf{x}^{\perp} \|^2 \right)
$, where $\bar{\mu}:=s_{\min}(\Psi)$. With $ \| \mathbf{x} ^{\|}\!-\! \mathbf{x'} \|^2 +   \| \mathbf{x}^{\perp} \|^2  =  \| \mathbf{x} \!-\! \mathbf{x'} \|^2$, this is  \eqref{eqboldRF_plus_cL_strgmon}. For $c>c_{\min}$ as in the statement, $\Psi \succ 0$,  $\bar{\mu}>0$. 
 $\hfill\blacksquare$ 

{\bf Proof of Lemma \ref{lem_monotone}: }
(i): The operator ${\bm{\mathfrak{A}}}$ \eqref{op_hat_A} is written as ${\bm{\mathfrak{A}}}\! =\! {\bm{\mathfrak{A}}}_1\!+\!{\bm{\mathfrak{A}}}_2$, where ${\bm{\mathfrak{A}}}_1\!\!=\!\!\mathcal{R}^TN_{\Omega} \mathcal{R}(\mathbf{x})\! \!\times  \!\!\bm{0}_{Nm}\! \!\times \! \!N_{\bm{R}^{Nm}_{+}}(\bm{\lambda})$ and 
${\bm{\mathfrak{A}}}_2$ the skew-symmetric matrix in \eqref{op_hat_A} (maximally monotone by \cite[Ex. 20.30]{combettes1}).  
$N_{\Omega}$ and $N_{\bm{R}^{Nm}_{+}}$ are maximally monotone (normal cones of closed convex sets, \cite[Thm. 20.40]{combettes1},  
as is $\bm{0}_{Nm}$. 
Since $\mathcal{R}$ is full row rank,  $\mathcal{R}^T \!\!\circ \!\!N_{\Omega}\! \!\circ \!\! \mathcal{R}$ is also maximally monotone, \cite[Prop.  12.5.5]{FacchineiBOOK}. 
As the direct sum of maximally monotone operators, ${\bm{\mathfrak{A}}}_1$ is maximally monotone, \cite[Prop.  20.23]{combettes1}, and since $dom {\bm{\mathfrak{A}}}_2\!=\! \bm{R}^{Nn+2Nm}$, ${\bm{\mathfrak{A}}}$ is also maximally monotone,  \cite[Cor. 24.4]{combettes1}.

(ii): Let $\varpi\!=\! col(\mathbf{x}, \mathbf{z}, \mathbf{\lambda})$, $\varpi' \!=\! col(\mathbf{x'}, \mathbf{z'}, \mathbf{\lambda'})$. For ${\bm{\mathfrak{B}}}$, \eqref{op_hat_A},
  \vspace{-0.2cm}
\begin{align}\label{eq_Bcocoercive_1}
& \hspace{-0.5cm}\!\langle \varpi -\varpi',\bm{\mathfrak{B}} \varpi  - \bm{\mathfrak{B}} \varpi' \rangle  \\
& \quad = \langle \mathbf{x}\!-\! \mathbf{x'},  \! \mathcal{R}^T\mathbf{F}(\mathbf{x})\!-\! \mathcal{R}^T\mathbf{F}(\mathbf{x'})\!+\!c \mathbf{L}_x \mathbf{x} \!-\! c \mathbf{L}_x \mathbf{x'} \rangle  \nonumber \\
& \quad \quad + \langle \bm{\lambda}- \bm{\lambda}',  \mathbf{L}_\lambda\, \bm{\lambda} - \mathbf{L}_\lambda \, \bm{\lambda}'   \rangle \nonumber 
\end{align}
Note that operator $\mathbf{L}_x$ is $\|\mathbf{L}_x\|$-Lipschitz, with    $ \|\mathbf{L}_x\| \!=\! s_N(L) \leq 2d^*$, 
Then, by the triangle inequality, using  
Lipschitz continuity of $\mathbf{F}$ (Assumption \ref{Lipchitz_boldFassump}) and $\|\mathcal{R}\| = 1$, we can write,   \vspace{-0.2cm}
\begin{eqnarray}\label{eqboldRF_plus_cL_Lip}
&&\hspace{-1cm}\! \| \mathcal{R}^T\mathbf{F}(\mathbf{x})\!-\!\mathcal{R}^T\mathbf{F}(\mathbf{x'})\!+\!c \mathbf{L}_x\mathbf{x} \!-\! c \mathbf{L}_x\mathbf{x'} \| \!\leq \! \bar{\theta}\|\mathbf{x}\!-\!\mathbf{x'}\|
\end{eqnarray}
for any $\mathbf{x}$ and $\mathbf{x}'$, where $\bar{\theta} \!= \! \theta \!+ \! 2 c  d^*$.  

By Lemma \ref{boldRF_plus_cL_strgmon_lemma}, for the first term on the right-hand side of \eqref{eq_Bcocoercive_1}, \eqref{eqboldRF_plus_cL_strgmon} holds for any $\mathbf{x}$ and any  $\mathbf{x'} \in \mathbf{E}_x$. 
Using this 
and \eqref{eqboldRF_plus_cL_Lip},  we can write, for any $\mathbf{x}$ and  any $\mathbf{x'} \in \mathbf{E}_x$,  \vspace{-0.2cm}
\begin{eqnarray}\label{eqboldRF_plus_cL_cocoercive}
&&\hspace{-1cm}\langle \mathbf{x}\!-\! \mathbf{x'},  \! \mathcal{R}^T\mathbf{F}(\mathbf{x})\!-\!\mathcal{R}^T\mathbf{F}(\mathbf{x'})\!+\!c \mathbf{L}_x \mathbf{x} \!-\! c \mathbf{L}_x \mathbf{x'} \rangle \nonumber \\
&& 
\!\geq  \frac{\bar{\mu}}{\bar{\theta}^2}  \| \mathcal{R}^T\mathbf{F}(\mathbf{x})\!+\!c \mathbf{L}_x\mathbf{x}\!-\!\mathcal{R}^T\mathbf{F}(\mathbf{x'})\!-\!c \mathbf{L}_x \mathbf{x'} \|^2.
\end{eqnarray}
For the second term on the right-hand side of \eqref{eq_Bcocoercive_1}, we note that since $\mathbf{L}_\lambda$ is symmetric,   $\mathbf{L}_\lambda  \, \bm{\lambda}$ 
 is the gradient of $\tilde{f}(\bm{\lambda})\!:=\!\frac{1}{2} \bm{\lambda}^T \! \mathbf{L}_\lambda \bm{\lambda}$, 
which is convex since $\nabla^2 \tilde{f}(\bm{\lambda})\!=\!\mathbf{L}_\lambda \succeq 0$  (by Assumption \ref{connectivity}).  
  As the $\|\mathbf{L}_\lambda\|-$Lipschitz continuous gradient of a convex function,  it follows that  $\mathbf{L}_\lambda  \bm{\lambda}$ is $\frac{1}{\|\mathbf{L}_\lambda \|}$-cocoercive, (by Baillon-Haddad theorem, \cite[Thm. 18.15]{combettes1}), 
where   $\frac{1}{\|\mathbf{L}_\lambda \|} \geq \frac{1}{2d^*}$. 
 Thus,  for  any $\bm{\lambda}$ and $\bm{\lambda}'$, \vspace{-0.2cm}  
  \begin{equation}\label{equ_lem5_3_2}
 \langle \bm{\lambda}- \bm{\lambda}',  \mathbf{L}_\lambda\, \bm{\lambda} - \mathbf{L}_\lambda \, \bm{\lambda}'   \rangle \geq 
 \frac{1}{2d^*} \| \mathbf{L}_\lambda \, \bm{\lambda} - \mathbf{L}_\lambda \, \bm{\lambda}' \|^2.
\end{equation} 
Then, for  any $\varpi$ 
and  any $\varpi' \in {\bm{\Omega}}_E$, since  $ \mathbf{x'} \in \mathbf{E}_x$, 
using \eqref{eqboldRF_plus_cL_cocoercive} and \eqref{equ_lem5_3_2} in \eqref{eq_Bcocoercive_1}, it follows that  \vspace{-0.2cm}
\begin{align*}
& \hspace{-0.5cm}\!\langle \varpi -\varpi',\bm{\mathfrak{B}} \varpi  - \bm{\mathfrak{B}} \varpi' \rangle  \\
& \quad = \langle \mathbf{x}\!-\! \mathbf{x'},  \! \mathcal{R}^T\mathbf{F}(\mathbf{x})\!-\! \mathcal{R}^T\mathbf{F}(\mathbf{x'})\!+\!c \mathbf{L}_x \mathbf{x} \!-\! c \mathbf{L}_x \mathbf{x'} \rangle  \nonumber \\
& \quad \quad + \langle \bm{\lambda}- \bm{\lambda}',  \mathbf{L}_\lambda\, \bm{\lambda} - \mathbf{L}_\lambda \, \bm{\lambda}'   \rangle \nonumber \\
& \quad \geq  \frac{\bar{\mu}}{\bar{\theta}^2}  \| \mathcal{R}^T\mathbf{F}(\mathbf{x})\!+\!c \mathbf{L}_x\mathbf{x}\!-\!\mathcal{R}^T\mathbf{F}(\mathbf{x'})\!-\!c \mathbf{L}_x \mathbf{x'} \|^2 \nonumber \\
&  \quad \quad  + \frac{1}{2d^*} \| \mathbf{L}_\lambda \, \bm{\lambda} - \mathbf{L}_\lambda \, \bm{\lambda}' \|^2 \nonumber \\
& \quad \geq   \min \{ \frac{\bar{\mu}}{\bar{\theta}^2},  \frac{1}{2d^*}  \} \| \bm{\mathfrak{B}} \varpi  - \bm{\mathfrak{B}} \varpi'\|^2
\end{align*}
 Hence, \eqref{eq:B_cocercive_restricted} holds for any $0<\!\beta \!\leq \!  \min \{ \frac{\bar{\mu}}{\bar{\theta}^2},  \frac{1}{2d^*}  \}$. 
\hfill $\Box$

{\bf Proof of Lemma \ref{lem_PhiA_B_prop}: }
(i) For any $(\varpi,u)\in gra {\Phi}^{-1}{\bm{\mathfrak{A}}}$ and $(\varpi',v)\in gra {\Phi}^{-1}{\bm{\mathfrak{A}}}$, ${\Phi} u \in {\Phi} {\Phi}^{-1}{\bm{\mathfrak{A}}}(\varpi)\in {\bm{\mathfrak{A}}}(\varpi) $ and ${\Phi} v \in {\Phi} {\Phi}^{-1}{\bm{\mathfrak{A}}}(\varpi')\in {\bm{\mathfrak{A}}}(\varpi') $. Since $\Phi \succ 0$ (cf. Lemma \ref{lem_monotone_metric}) and ${\bm{\mathfrak{A}}}$ is maximally monotone (cf. Lemma \ref{lem_monotone}(i)), $\langle \varpi-\varpi',u-v\rangle_{{\Phi}}= \langle  \varpi-\varpi', {\Phi}(u-v) \rangle \geq 0, \forall \varpi,\varpi'\in dom \mathfrak{A}$, and ${\Phi}^{-1}{\bm{\mathfrak{A}}}$ is maximally monotone  under the ${\Phi}-$induced inner product. 
By  \cite[Prop.  23.7]{combettes1},  $T_2= ({\rm Id}+\Phi^{-1}\mathfrak{A})^{-1}$ is firmly nonexpansive under  the $\Phi-$induced norm $\| \cdot\|_{\Phi}$,  hence $T_2 \in\mathcal{A}(\frac{1}{2}) $. \\
(ii) By Lemma \ref{lem_monotone_metric},  ${\Phi} - \delta I_{n+2Nm} \succeq 0$, hence 
 the eigenvalues of ${\Phi}$ satisfy  $s_{\max}({\Phi})\geq s_{\min}({\Phi}) \geq \delta$.
By  \eqref{eq:B_cocercive_restricted} in Lemma \ref{lem_monotone},  for any  $\varpi $ and  any $\varpi' \in \bm{\Omega}_E$ it holds that, 
 \vspace{-0.2cm}
\begin{equation}\label{eq:Bcocercive_Phi_0}
\begin{array}{ll}
\hspace{-0.3cm} \langle  \varpi\!-\!\varpi' ,\! \Phi^{-1}{\bm{\mathfrak{B}}}\varpi\!-\!\Phi^{-1}{\bm{\mathfrak{B}}}\varpi' \rangle_{\Phi} \!&\!\! \!=\! \langle {\bm{\mathfrak{B}}}\varpi\!-\!{\bm{\mathfrak{B}}}\varpi',  \!\varpi\!-\!\varpi'  \rangle \\
& \!\!\!\geq \! \beta  \|{\bm{\mathfrak{B}}}\varpi\!-\!{\bm{\mathfrak{B}}}\varpi' \|^2
\end{array}
\end{equation} 
Since ${\bm{\mathfrak{B}}}$ is single-valued and ${\Phi}\succ0$,   \vspace{-0.2cm} 
\begin{equation}
\begin{array}{l}\nonumber
s_{\max}(\Phi^{-1})  \|{\bm{\mathfrak{B}}}\varpi\!-\!{\bm{\mathfrak{B}}}\varpi' \|^2 
\geq \! \langle \!{\bm{\mathfrak{B}}}\varpi \!- \!{\bm{\mathfrak{B}}}\varpi',{\Phi}^{-1}({\bm{\mathfrak{B}}}\varpi \!- \!{\bm{\mathfrak{B}}}\varpi') \! \rangle \\
\qquad = \langle {\Phi}{\Phi}^{-1}({\bm{\mathfrak{B}}}\varpi-{\bm{\mathfrak{B}}}\varpi'),{\Phi}^{-1}({\bm{\mathfrak{B}}}\varpi-{\bm{\mathfrak{B}}}\varpi') \rangle \\
\qquad = \| {\Phi}^{-1}{\bm{\mathfrak{B}}}\varpi-{\Phi}^{-1}{\bm{\mathfrak{B}}}\varpi' \|^2_{{\Phi}}.
\end{array}
\end{equation}
Using $s_{\max}(\Phi^{-1}) = 1/s_{\min}(\Phi)  \leq 1/\delta$, it follows that $
 \|{\bm{\mathfrak{B}}}\varpi\!-\!{\bm{\mathfrak{B}}}\varpi' \|^2 \geq \delta\| {\Phi}^{-1}{\bm{\mathfrak{B}}}\varpi-{\Phi}^{-1}{\bm{\mathfrak{B}}}\varpi' \|^2_{{\Phi}}$. 
 Using this in \eqref{eq:Bcocercive_Phi_0}, yields that, for any  $\varpi $ and  any {$\varpi' \in \bm{\Omega}_E$,  
$\langle \! \varpi\!-\!\varpi' ,\!  \Phi^{-1}{\bm{\mathfrak{B}}}\varpi\!-\!\Phi^{-1}{\bm{\mathfrak{B}}}\varpi'  \! \rangle_{\Phi} 
 \!\geq \! {\beta}{\delta}  \| \Phi^{-1}{\bm{\mathfrak{B}}}\varpi\!-\!\Phi^{-1}{\bm{\mathfrak{B}}}\varpi' \|^2_{\Phi},
$ hence  $\Phi^{-1}{\bm{\mathfrak{B}}}$ is $\beta\delta$-\textit{restricted} cocoercive. 
Then, for $T_1=Id \! -\!\Phi^{-1}\bm{\mathfrak{B}}$, using this property,  we can write  for any $ \varpi $ and any $\varpi' \in {\bm{\Omega}}_E$, \vspace{-0.2cm} 
\begin{equation*}
\begin{array}{ll}
 \hspace{-0.2cm}\| T_1 \varpi\! -  \!T_1 \varpi' \|^2_{\Phi} \!\!&= \| \varpi - \varpi' -\Phi^{-1}{\bm{\mathfrak{B}}} \varpi + \Phi^{-1}\bm{\mathfrak{B}} \varpi' \|^2_{\Phi}  \\
& =  \| \varpi - \varpi'\|^2_{\Phi} + \|\Phi^{-1}{\bm{\mathfrak{B}}} \varpi - \Phi^{-1}\bm{\mathfrak{B}} \varpi' \|^2_{\Phi} \\
&\quad - 2\langle \varpi - \varpi', \Phi^{-1}{\bm{\mathfrak{B}}} \varpi - \Phi^{-1}\bm{\mathfrak{B}} \varpi' \rangle_\Phi \\
& \leq \| \varpi - \varpi'\|^2_{\Phi} + \|\Phi^{-1}{\bm{\mathfrak{B}}} \varpi - \Phi^{-1}\bm{\mathfrak{B}} \varpi' \|^2_{\Phi} \\
&   \quad                           - 2 \beta \delta \|\Phi^{-1}{\bm{\mathfrak{B}}} \varpi - \Phi^{-1}\bm{\mathfrak{B}} \varpi' \|^2_{\Phi} \\
\end{array}
\end{equation*}
which is \eqref{eq:prop1_T1}. Hence, 
since $2\beta\delta >1$ by assumption,  $\| T_1 \varpi -  T_1 \varpi' \|^2_{\Phi} \leq  \| \varpi - \varpi' \|^2_{\Phi}$ ($T_1$ is \textit{restricted} nonexpansive).
\hfill $\Box$


\begin{IEEEbiography}[{\includegraphics[width=0.9in]{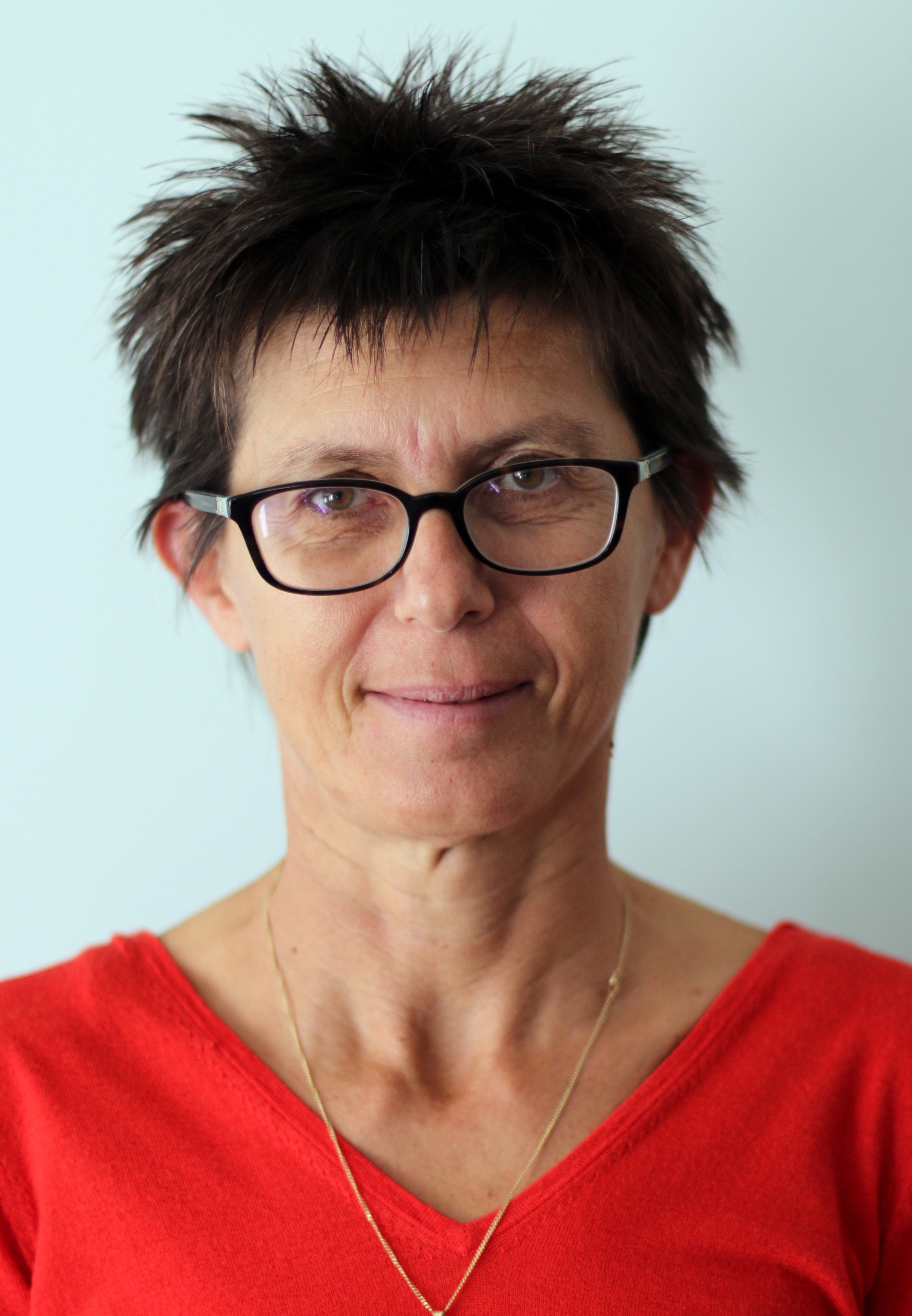}}]{Lacra Pavel}
(M'92 - SM'04)  received the Diploma of Engineer from Technical University of Iasi, Romania and the Ph.D. degree in Electrical Engineering from Queen's University at Kingston, Canada.  After a postdoctoral stage at the National Research Council and four years of working in the industry, she joined University of Toronto, Canada in August 2002, where she is a Professor in the Department of Electrical and Computer Engineering. Her research interests are in game theory and distributed optimization in networks, with emphasis on dynamics and control aspects. She acted as Publications Chair of the 45th IEEE Conference on Decision and Control. She is the author of  the book {\em Game Theory for Control of Optical Networks} (Birkh\"{a}user-Springer Science, ISBN 978-0-8176-8321-4, 2012).
\end{IEEEbiography}

\end{document}